\documentclass[titlepage]{article} % This command is used to set the type of document you are working on such as an article, book, or presenation

\usepackage{geometry} % This package allows the editing of the page layout
\usepackage{amsmath}  % This package allows the use of a large range of mathematical formula, commands, and symbols
\usepackage{graphicx}  % This package allows the importing of images
\usepackage{amssymb}
\usepackage{cite}
\usepackage{float}
\usepackage{tikz}
\usetikzlibrary{calc}
\usepackage{tkz-tab}
\usepackage{caption}
\usepackage{latexsym}
\usepackage{subcaption}
\usepackage{hyperref}
\hypersetup{linkbordercolor=blue}
\usetikzlibrary{decorations.pathreplacing}

\newcommand{\maketitletwo}[2][]{\begin{center}
	\Large{\textbf{Differentiable Lagrangian Shock Hydrodynamics with application to stable shock acceleration of density interfaces}
	} % Name of course here
	\vspace{5pt}

	\normalsize{Kevin Korner*, Brandon Talamini, Julian Andrej, Michael Tupek, Bill Moses, \\ Rob Rieben, Tzanio Kolev, Jamie Bramwell, Dan White, Jon Belof, William Schill \\} % % Your name here
	\normalsize{*Corresponding author: korner1@llnl.gov \\}
	\normalsize{Lawrence Livermore National Laboratory\\}
	\normalsize{LLNL-JRNL-872441}

\end{center}}
\begin{document}
% \maketitle
\maketitletwo[5]

% \tableofcontents
% \newpage

\begin{abstract}

We develop a gradient based optimization approach for the equations of compressible, Lagrangian hydrodynamics and
demonstrate how it can be employed to automatically uncover strategies to control hydrodynamic instabilities arising from
shock acceleration of density interfaces. Strategies for controlling the Richtmyer-Meshkov instability (RMI) are of great benefit
for inertial confinement fusion (ICF) where shock interactions with many small imperfections in the density interface 
lead to instabilities which rapidly grow over time. These instabilities lead to mixing which, in the case of laser driven ICF, quenches the runaway fusion process ruining the potential for positive energy return.
We demonstrate that control of these instabilities can be achieved by optimization of initial conditions with ($>100$) parameters. Optimizing over a large parameter space like this is not possible with gradient-free optimization strategies.
%judicious
This requires computation of the gradient of the outputs of a numerical solution to the equations of Lagrangian hydrodynamics with respect to the inputs.
We show that the efficient computation of these gradients is made possible via a judicious application of (i) adjoint methods, the exact formal representation of sensitivities involving partial differential equations, and (ii) automatic differentiation (AD), the algorithmic calculation of derivatives of functions.
Careful regularization of multiple operators including artificial viscosity and timestep control is required.
We perform design optimization of $>100$ parameter energy field driving the Richtmyer Meshkov instability showing significant suppression while simultaneously enhancing the acceleration of the interface relative to a nominal baseline case. 

\end{abstract}

\section{Introduction}

Shock hydrodynamics describe the behavior of some of the most complex phenomena in science and technology. 
Examples abound including supernovae in astrophysics and inertial or magnetic confinement fusion in the laboratory. 
The physical processes that may occur are manifold; however, much of the complexity may be thought of as arising from dynamical instabilities. 
By way of example, the Richtmyer-Meshkov instability (RMI) occurs when a shock wave -- a discontinuous jump in the thermodynamic state of the material -- impinges on an interface between two materials of differing densities\cite{richtmyer1954taylor,taylor1950instability}.
The resulting dynamics are famously unconditionally unstable resulting in \textit{jetting} which rips up the interface and mixes the materials together.
For a comprehensive review of RMI, see~\cite{zhou_rayleightaylor_2017,zhou_rayleightaylor_2017-1} and the references therein.
RMI has held a critical role in both scientific and technological applications including astrophysics~\cite{zhou2021rayleigh}, mining~\cite{birkhoff1948explosives}, many applications of fluid transport~\cite{zhou2019turbulent} including scramjets~\cite{zhou2021rayleigh}, and laser driven inertial confinement fusion (ICF)~\cite{mikaelian2011extended,remington2019rayleigh} such as is pursued at the National Ignition Facility (NIF). 
Control over RMI induced jetting is thus a grand challenge.

Significant recent effort has gone into optimization involving hydrodynamic control over jetting including gas gun driven experiments~\cite{sterbentz2022design}, and explosively driven linear shaped charges~\cite{kline2024reducing,sterbentz2024explosively}.
These strategies typically proceed by identifying a \textit{small} number of parameters describing a physical system of interest, performing a large suite of calculations, fitting a machine-learned model to predict certain interesting quantities from a simulation, and then using the resulting model to optimize or otherwise interrogate the behavior of the system~\cite{jekel2024machine,sterbentz2023linear}. 
It is even possible to develop an analytic solution~\cite{schill_suppression_2024} to a reduced and simplified version of this problem.
A closely related class of problems is the inference of material parameters from complex extreme pressure materials science experiments such as those conducted on pulsed power~\cite{schill2021simultaneous,schill2023inference} or laser~\cite{gorman2023ramp} platforms.
These techniques are powerful but they are currently limited to a small number of design variables; however, the general case -- i.e. control over an arbitrary interface with potentially $>100$ degrees of freedom remains unsolved.
The fundamental reason for this has to do with basic characteristics of optimization theory; optimization is very expensive for large numbers of variables in the absence of gradient information and the computation of the objective function gradients in these problems rely on the \textit{derivative of the hydrocode outputs with respect to their inputs}.
Computing this gradient in the context of Lagrangian Shock Hydrodynamics is a challenging task and is one of the key contributions of this paper.

The subject of gradient based optimization involving partial differential equations has seen extensive study.  For instance, topology optimization is a field of engineering and mathematics which deals with optimizing the material layout within a given design space to achieve the best performance under specified conditions and has been applied to studying many different problems~\cite{akerson_optimal_2022, akerson_optimal_2023, bendsoe_topology_2004}. 
Basic adjoint calculations can be found in a variety of texts (see, for instance,~\cite{akerson_optimal_2023,plessix_review_2006}). 
This has seen significant recent development by some of the authors~\cite{doecode_94261} for nonlinear problems in solid mechanics.
Recently, a number of authors have begun applying these techniques to challenging problems in high energy physics~\cite{carli2023algorithmic}.

The remainder of this article will be organized as follows: (i) gradient based optimization of time dependent problems (Section~\ref{sec:TimeDependentGradient}), (ii) discretization by finite elements of Lagrangian shock hydrodynamics (Section~\ref{sec:FiniteElements}), and (iii) the application of this technique to the suppression of RMI suggesting a tantalizing pathway to stable shock acceleration of density interfaces (Section~\ref{sec:LagrangianHydrodynamics}).

\section{Gradient based optimization of time dependent problems}\label{sec:TimeDependentGradient} 

In this section, we introduce several concepts that are necessary for our ultimate goal
of computing gradients of a Lagrangian hydrodynamics discretization. 
First, we discuss the two main strategies used to develop adjoints of time dependent problems: discretize the differentiation and differentiate the discretization.
We argue for the latter case. 
Next we discuss methods for verifying the computation of gradients by  checking the error order of convergence. 
We then introduce automatic differentiation (AD), quantities of interest (QOIs) and check-pointing for time dependent problems. 
Finally, we consider a simple example involving multi-particle systems to illustrate the combination of these components.

%We begin by covering two methods of calculating the adjoints each having advantages and disadvantages. 
%Another way to describe the difference is in how we treat the time propagation. 

\subsection{Discretizing the Differentiation}\label{sec:TraditionalAdjoint}
One method of conducting an adjoint calculation is to \textbf{\textit{discretize the differentiation}}~\cite{plessix_review_2006}. 
We find a system of equations which can be solved to find the derivative, then we discretize the system in time in order to do calculations. 
Due to the fact that no discretization was necessary in order to describe both the forward and adjoint problems mathematically, it presents itself nicely for theoretical work in the subject. 
This makes it ideal for studying properties of solutions using analytical methods. 
Unfortunately, the vast majority of problems are not analytically integrable in time, so, in order to represent the solution, we must choose time integrators for both the forward and reverse problems. 
A problem arises because, in principle, if the forward and backward time integrators are not chosen specifically to be exact adjoints, a small error will be present at every adjoint solve. 
This error will then accrue as we propagate the adjoint and become non-negligible in the final derivative. 

Consider a dynamical system given by 
\begin{equation}
	\dot{y} = f(y)	
	\label{eq:dynamicalsystem}
\end{equation}
where $y \in \mathbb{R}^N$, $f : \mathbb{R}^N \to \mathbb{R}^N$, and $\dot{()}$ indicates a time derivative.

We define a generic optimization problem as 
\begin{align*}
	O[y_0] &= \int_0^T o(y(t), t) dt \, \\
	y(0) &= y_0 \, \\	
	\dot{y} &= f(y)\, ,
\end{align*}
where $o$ is a time dependent objective function, which is a measure of quality of the solution, and $y_0$ is the solution at the initial time step. We would like to minimize the objective function with respect to the initial state $y_0$. The difficulty here arises from the fact that the integrated solution must satisfy the dynamics. In other words, the quantity
\begin{equation}
	\dfrac{\partial y}{\partial y_0}(t)
\end{equation}
is not trivial to calculate. As a result, the adjoint method was developed. The first step is to define an augmented objective function of the form
\begin{equation}
	\tilde{O}[y_0, \lambda] = \int_0^T o(y(t), t) - \lambda(t) \cdot (\dot{y}(t) - f(y(t))) dt
\end{equation}
where $\lambda(t)$ is known as the adjoint variable. Now, if we perturb our objective with respect to the variables of interest, we have
\begin{equation}
	\delta \tilde{O}(y_0, \lambda) = \int_0^T \left( \dfrac{\partial o}{\partial y} \cdot \left( \dfrac{\partial y}{\partial y_0} \cdot \delta y_0 \right) - \lambda(t) \cdot \left( \dfrac{d}{dt} \left(\dfrac{d y}{dy_0}\delta y_0 \right)  - \dfrac{\partial f}{\partial y} \dfrac{\partial y}{\partial y_0}\delta y_0\right)\right)dt\, .
	\label{eq:augmentedlagrangian2}
\end{equation}

For simplicity of notation, we write $\dfrac{\partial y}{\partial y_0} \cdot \delta y_0 = y_{y_0}$. Rearranging and applying integration by parts on the time derivative, we have
\begin{align*}
	\delta \tilde{O}(y_0, \lambda) = \int_0^T \left( \dfrac{\partial o}{\partial y} + \dot{\lambda} + \lambda \cdot \dfrac{\partial f}{\partial y} \right)\cdot y_{y_0} dt - \left( \lambda(t)\cdot \dfrac{\partial y}{\partial y_0} \delta y_0 \right)_{T = 0}^{T=T}\, .
\end{align*}

Now, we make a careful choice of $\lambda$, such that it satisfies
\begin{align}
	\dot{\lambda} &= - \lambda \cdot \dfrac{\partial f}{\partial y} - \dfrac{\partial o}{\partial y} \\
	\lambda(T) &= 0\, .
	\label{eq:continuoussensitivity}
\end{align}

Then, Equation~(\ref{eq:augmentedlagrangian2}) simplifies to

\begin{equation}
	\delta \tilde{O}(y_0, \lambda) = \lambda(0) \cdot \delta y_0 \quad \to \quad \dfrac{\partial \tilde{O}}{\partial y_0} = \lambda(0)\, .
\end{equation}

\subsection{Differentiating the Discretization}\label{sec:graphnetworkapproach}
%The other method (graph network approach), operates by taking the derivative of the discretized system directly. 
%
In this work, we argue its best to \textbf{\textit{differentiate the discretization}} for the case of non-linear partial
differential equations representing temporal evolution of conserved physical quantities. 
Compared to the method presented in the previous section, nothing changes about the forward problem. 
The same discretization errors that exist in that context exist in this. 
The real benefit comes in that there is no choice in the integrator for the reverse pass. 
The solution is given by the adjoint of the time integrator used by the forward problem. 
Because we are directly differentiating the time stepping algorithm used for the forward problem, we are technically differentiating the computational graph directly.
We will refer to this process as the graph network approach, because we operate on the computational graph rather than the continuous dynamics.
For the case of Lagrangian hydrodynamics, there are important aspects of the forward time integration process which need to be preserved like exact conservation of total
energy. 
A possible interpretation of this is that the graph network approach calculates the optimal backwards time stepping scheme. 
The main drawback of this approach is that it is harder to do rigorous proofs of convergence about the continuous problem because we operate on the discretized domain.

The graph network approach builds off of the time discretized system. We take the dynamical system given by Equation~(\ref{eq:dynamicalsystem}). Then, composing the time derivative equation with an update rule gives a discrete time update operator of the form
\begin{equation}\label{eq:discretedynamics}
	y_{i+1} = G(y_i, dt)\, ,	
\end{equation}
where $G$ denotes a particular form of time integrator such as Euler timestepping or energy conserving two-stage Runge-Kutta. 
Now, when we
formulate the optimization problem, we can take a much more direct approach. We now write the time integral as a discrete sum over each time step as
\begin{align*}
	O[y_0] &= \sum_{k = 0}^{N_t - 1} o(y_k ,t_k) \Delta t_k
\end{align*}
Now, we take the derivative with respect to the initial configuration $y_0$ as 
\begin{align*}
	\dfrac{\partial O}{\partial y_0} = \sum_{k = 0}^{N_t - 1} \dfrac{\partial o}{\partial y_k}(y_k, t_k) \cdot \dfrac{\partial y_k}{\partial y_0}\Delta t_k\, .
\end{align*}
There are a few properties that we can note from looking at this equation. First of all, if we start from the right most quantities and apply the 
operations %matrices%% matrices? we havent discretized yet 
iteratively, we end up with a problem of poorly controlled growth of computational cost. 
Many discrete dynamical systems represent these operations as matrices (particularly those in finite elements) which have a sparse structure arising from locality of shape functions used for the discretization. 
This property is something we would like to take advantage of when finding solutions. 
When multiplying sparse matrices together, the resultant matrix is often less sparse than either of the two original matrices. 
After a few iterations of sparse matrix compositions, we will then have a dense matrix which we would not like to store in memory.
Many implementations of finite element systems also prefer to be matrix free (matrices only exist in their action on an applied vector). The sums, however, show a very self-similar structure which can be manipulated to make the above calculation both fast and efficient.
Note that we can rearrange the sum as
\begin{align*}
	\dfrac{\partial O}{\partial y_0} = \dfrac{\partial o}{\partial y_{N_t - 1}} \cdot \dfrac{\partial y_{N_t - 1}}{\partial y_{N_t - 2}} \cdot \dfrac{\partial y_{N_t - 2}}{\partial y_{N_t - 3}} \cdots \dfrac{\partial y_2}{\partial y_1} \cdot \dfrac{\partial y_1}{\partial y_0} \Delta t_{N_t - 1}&\\	
	+ \dfrac{\partial o}{\partial y_{N_t - 2}} \cdot \dfrac{\partial y_{N_t - 2}}{\partial y_{N_t - 3}} \cdots \dfrac{\partial y_2}{\partial y_1} \cdot \dfrac{\partial y_1}{\partial y_0}\Delta t_{N_t - 2}&\\
	+\quad \quad  \vdots \quad \quad &\\
	+ \dfrac{\partial o}{\partial y_2} \cdot \dfrac{\partial y_2}{\partial y_1} \cdot \dfrac{\partial y_1}{\partial y_0} \Delta t_{2}&\\
	+ \dfrac{\partial o}{\partial y_1} \cdot \dfrac{\partial y_1}{\partial y_0} \Delta t_1&\\
	+ \dfrac{\partial o}{\partial y_0}\Delta t_0&
\end{align*}
From here, we see we can construct the iteration
\begin{align}
	\lambda_{N_{t - 1}} &= \dfrac{\partial o}{\partial y_{N_t - 1}} \Delta t_{N_t - 1} \\
	\lambda_{k} &= \lambda_{k+1} \cdot \dfrac{\partial y_{k+1}}{\partial y_k} + \dfrac{\partial o}{\partial y_k}\Delta t_k	
\end{align}
This method works from the final time step and, much like the traditional adjoint method described in Section~(\ref{sec:TraditionalAdjoint}), we accumulate the derivative by stepping backwards in time. 
Additionally, by propagating the solution backwards in time, we only need to consider matrix-vector products, storing vectors, leading to fast and efficient calculations. A caveat is that every state of the system throughout the time series must be stored and accessible after the forward solve. 
We address this issue by using a checkpointing algorithm to optimally reconstruct states. 
The methods used will be described in Section~\ref{sec:checkpointing}. 
What remains is the ability to calculate the quantity
\begin{align*}
	\lambda \cdot \dfrac{\partial y_{k+1}}{\partial y_{k}}	 \, .
\end{align*}

%\subsection{Final Timestep Only}
A particular subset of problems of interest is when the objective function is only a function of the final timestep, written as
\begin{align*}
	O = O(y_f)\, .
\end{align*}
As a result, the adjoint structure also simplifies as 
\begin{align*}
	\lambda_{N_{t-1}} &= \dfrac{\partial O}{\partial y_{N_{t-1}}} \, , \\ 
	\lambda_k &= \lambda_{k+1} \cdot \dfrac{\partial y_{k+1}}{\partial y_k}
\end{align*}
Note that this can be done by utilizing the previous framework and setting
\begin{align*}
	o(y_k, t_k) = \delta_{k,N_{t-1}} \frac{O(y_k)}{\Delta t_k}
\end{align*}
Which type of objective function to use (or a combination of both) will depend on the nature of the problem to be solved.

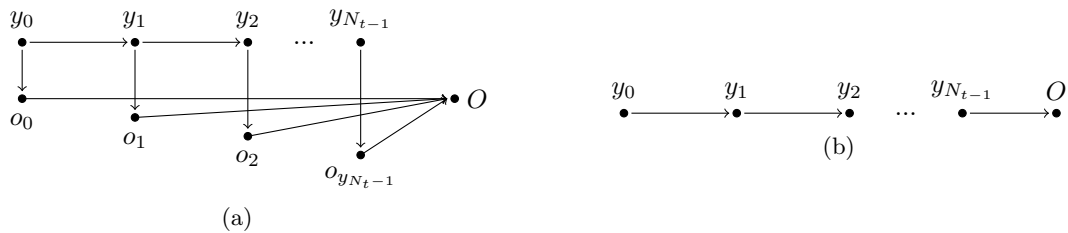
\begin{figure}[!hbt]
	\centering
	\begin{tikzpicture}
		\draw(0, 0) node[inner sep=0]{\begin{subfigure}[b]{.4\textwidth}
			\centering
			\begin{tikzpicture}
				\coordinate (A) at (0,0);
				\coordinate (B) at (1.5,0);
				\coordinate (C) at (3.0,0);
				\coordinate (D) at (4.5,0);
				\coordinate (o0) at (0,-.75);
				\coordinate (o1) at (1.5,-1.);
				\coordinate (o2) at (3.0,-1.25);
				\coordinate (oyn) at (4.5,-1.5);
				\coordinate (O) at (5.75,-0.75);
				% \coordinate (E) at (6.0,0);
				\draw[->] ($(A) + (0.1, 0)$) -- ($(B) - (0.1, 0)$);
				\draw[->] ($(B) + (0.1, 0)$) -- ($(C) - (0.1, 0)$);
				\draw[->] ($(A) + (0.0, -0.1)$) -- ($(o0) + (0.0, 0.1)$);
				\draw[->] ($(B) + (0.0, -0.1)$) -- ($(o1) + (0.0, 0.1)$);
				\draw[->] ($(C) + (0.0, -0.1)$) -- ($(o2) + (0.0, 0.1)$);
				\draw[->] ($(D) + (0.0, -0.1)$) -- ($(oyn) + (0.0, 0.1)$);

				\draw ($(A) + (0.0, 0.3)$) node {$y_0$};
				\draw [black,fill=black](A) circle(0.05); 
				\draw ($(B) + (0.0, 0.3)$) node {$y_1$};
				\draw [black,fill=black](B) circle(0.05); 
				\draw ($(C) + (0.0, 0.3)$) node {$y_2$};
				\draw [black,fill=black](C) circle(0.05); 

				\draw ($(D) + (0.0, 0.3)$) node {$y_{N_{t-1}}$};
				\draw [black,fill=black](D) circle(0.05);
				\draw ($(o0) - (0, 0.3)$) node {$o_0$};
				\draw [black,fill=black](o0) circle(0.05); 
				\draw ($(o1) - (0, 0.3)$) node {$o_1$};
				\draw [black,fill=black](o1) circle(0.05); 
				\draw ($(o2) - (0, 0.3)$) node {$o_2$};
				\draw [black,fill=black](o2) circle(0.05); 
				\draw ($(oyn) - (0, 0.3)$) node {$o_{y_{N_t - 1}}$};
				\draw [black,fill=black](oyn) circle(0.05); 
				% \draw ($(E) - (0, 0.3)$) node {$G$};
				% \draw [black,fill=black](E) circle(0.05); 
				\draw (3.75, 0) node {...};

				\draw ($(O) + (0.3, 0)$) node {$O$};
				\draw [black, fill=black](O) circle(0.05);
				\draw[->] ($(o0)$) -- ($(O) + (-0.1, 0.0)$);
				\draw[->] ($(o1)$) -- ($(O) + (-0.1, 0.0)$);
				\draw[->] ($(o2)$) -- ($(O) + (-0.1, 0.0)$);
				\draw[->] ($(oyn)$) -- ($(O) + (-0.1, 0.0)$);
			\end{tikzpicture}
			\caption{}
		\end{subfigure}
		};	
		\draw(8, 0) node[inner sep=0]{\begin{subfigure}[b]{.4\textwidth}
			\centering
			\begin{tikzpicture}
				\coordinate (A) at (0,0);
				\coordinate (B) at (1.5,0);
				\coordinate (C) at (3.0,0);
				\coordinate (D) at (4.5,0);
				\coordinate (O) at (5.75,0.0);
				% \coordinate (E) at (6.0,0);
				\draw[->] ($(A) + (0.1, 0)$) -- ($(B) - (0.1, 0)$);
				\draw[->] ($(B) + (0.1, 0)$) -- ($(C) - (0.1, 0)$);
				\draw[->] ($(D) + (0.1, 0)$) -- ($(O) - (0.1, 0)$);

				\draw ($(A) + (0.0, 0.3)$) node {$y_0$};
				\draw [black,fill=black](A) circle(0.05); 
				\draw ($(B) + (0.0, 0.3)$) node {$y_1$};
				\draw [black,fill=black](B) circle(0.05); 
				\draw ($(C) + (0.0, 0.3)$) node {$y_2$};
				\draw [black,fill=black](C) circle(0.05); 

				\draw ($(D) + (0.0, 0.3)$) node {$y_{N_{t-1}}$};
				\draw [black,fill=black](D) circle(0.05);
				\draw (3.75, 0) node {...};

				\draw ($(O) + (0.0, 0.3)$) node {$O$};
				\draw [black, fill=black](O) circle(0.05);
			\end{tikzpicture}
			\caption{}
		\end{subfigure}
		};	
	\end{tikzpicture}
	\caption{Graph networks for the forward solve of the time dependent objective (a) and the final timestep objective (b).}
\end{figure}

\subsection{Assessment of gradient behavior}

Our metric for what we consider to be a ``good'' gradient vs a ``bad'' one is whether it satisfies the famous Taylor remainder convergence test. Assume we are given a function $f$ and another function $df$ which is claimed to be the derivative of $f$. We can verify this claim by Taylor expanding it and checking the convergence. 
\begin{align}
	f(x + h \delta x) = f(x) + h df(x)\cdot \delta x + O(h^2) \, ,	
\end{align}
if we choose $\delta x$ to be normalized and $h \in \mathbb{R}$ to be small. This can also be written as
\begin{align}\label{eq:taylor_test}
	T(x, h, \delta x) = \|f(x + h \delta x) - f(x) - h df(x) \cdot \delta x \| = O(h^2)\, .
\end{align}
Therefore, in order to verify whether the given function $df$ is the derivative of the function $f$, we verify not only that the error goes to zero, but that we have quadratic convergence of the left hand side of Equation~\ref{eq:taylor_test}. Additionally, we verify that the direct error in the gradient ($h T(x, h, \delta x)$) also goes to zero.

\begin{figure}[!hbt]
	\centering	
	\includegraphics[width=.5\textwidth]{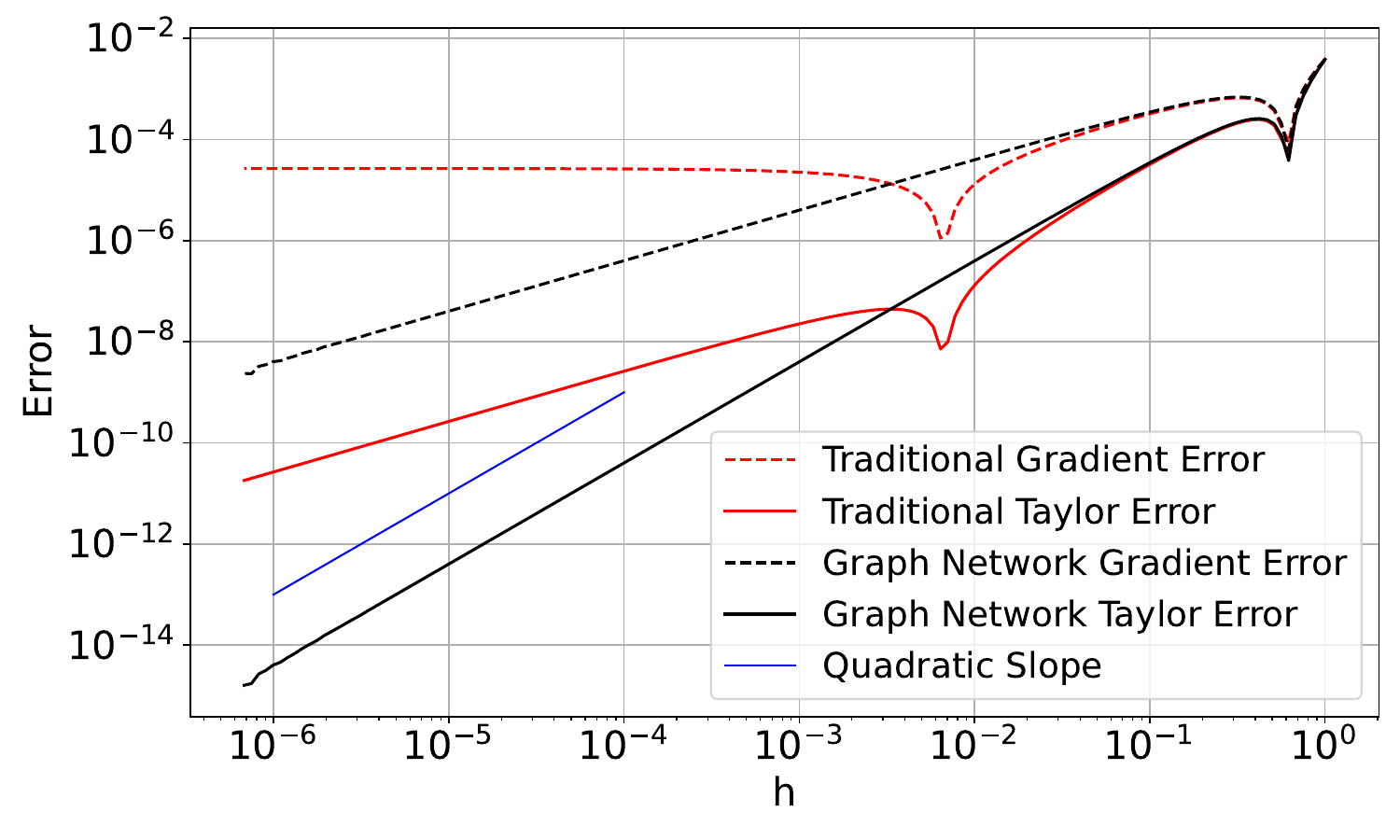}
	\caption{Error for the Taylor test for various perturbation magnitudes $h$. The black lines are the Taylor Error while the red are for the gradient error. The blue indicates expected scaling (quadratic) for the Taylor Error to confirm correct derivatives.}\label{fig:taylor_error}
\end{figure}

% \begin{align}
% 	\dot{\lambda} &= - \lambda \cdot \dfrac{\partial f}{\partial y} \\
% 	\frac{\lambda_{i} - \lambda_{i - 1} }{\Delta t} &= - \lambda_{i} \cdot \dfrac{\partial f}{\partial y}(y_i) \\ 
% 	-\lambda_{i-1} &= - \lambda_i - \Delta t \lambda_i \cdot \dfrac{\partial f}{\partial y} \\
% 	\lambda_{i-1} &= \lambda_i + \Delta t \lambda_i \cdot \dfrac{\partial f}{\partial y}
% \end{align}
To demonstrate, consider a nonlinear spring with dynamics given by $\ddot{y}(t) + k y(t) + \alpha y(t)^3 + b \dot{y}(t) = 0$. 
The parameters $k$, $\alpha$, and $\beta$ are the spring constant, hardening, and viscosity parameters, respectively. 
For simplicity, we set $k = \alpha = b = 1$. To integrate the system in time, we use a forward Euler integration scheme. 
For the backwards solve of the traditional adjoint, we use a backwards Euler scheme which uses the forward step. 
We integrate this system with initial conditions of $y(0) = 1.1$, $\dot{y}(0) = 5.0$ and integrate until $T = 10$, define a scalar objective function $O(y_f) = y(T) + \dot{y}(T)$, and conduct the above Taylor test and plot both the gradient error and the Taylor error in Figure~\ref{fig:taylor_error}.

As can be seen, the Graph Network approach has both smaller gradient error, and the proper convergence compared to the Traditional method. 
Additionally, the Traditional method gradient doesn't properly converge in the correct sense. 
This can present many problems in higher order optimization schemes, as accurate gradients are necessary.

Note that the discrepancy came in the choice of backwards integrator for the reverse problem in the Traditional method. 
While choosing the ``correct" integrator when using a forward Euler method can easily be remedied, when using more complex time integrators, namely higher order Runge-Kutta schemes, the choice is not as obvious. 
The graph network approach, consequently, presents a significant advantage over the traditional approach in that the choice of integrator for the reverse pass is exactly specified.

\subsection{Automatic Differentiation}
Automatic differentiation (also known as ``auto-diff'', ``auto-grad'', or simply ``AD'') is a computational technique used to efficiently and accurately evaluate derivatives of mathematical functions. 
It plays a crucial role in various fields such as machine learning, optimization, physics, simulations, and scientific computing.
The most popular AD libraries are PyTorch\cite{paszke2017automatic}, TensorFlow\cite{tensorflow2015-whitepaper}, and Jax\cite{jax2018github}, specifically due to their integrations into popular machine learning libraries.

The basic idea behind AD is to decompose functions into a sequence of elementary operations (such as addition, multiplication, exponentiation, etc.), for which the derivatives are well known.
Then by applying the chain rule iteratively to these operations, AD can compute the derivative of the entire function with respect to its input variables.
This technology enables rapid prototyping of tools involving differentiation, and can readily differentiate complex functions.

AD is extremely useful when developing analytical tools. 
Consider the case where AD is currently being used extensively: machine learning. 
Without AD tools, users of ML packages such as PyTorch and Tensorflow would have to manually specify gradients of their loss functions with respect to the neural network architecture. 
While this is not an impossible task, it would put severe limitations on who can use these tools, how quickly different models can be tested, and introduce many avenues for error. 
Arguably, the integration of automatic differentiation is the impetus which allowed the field to thrive as it has today.

AD is also important in expanding the scope of problems we can study. 
In many cases, especially problems with iteration, composition, recursion, branches, and complex algebras, taking derivatives by hand is not feasible and, as with the previous point, extremely prone to error. 
A much more convenient approach is to define computational methods which can handle the complexity and accurately traverse the computational graphs.

Automatic differentiation, however, is not without its faults and drawbacks. 
Often, naive implementations of AD can lead to inefficient run times and massive memory usages. 
This is because, unless otherwise specified, the tool must hold the entire computational graph, including sensitivities, in memory all at once. 
In large scale problems, particularly those involving dynamic simulations of hydrodynamic systems, this is cost prohibitive. 
As a result, we must control how we apply automatic differentiation and carefully consider various aspects such as memory allocation and code structure.

For the following work, we use the automatic differentiation library Enzyme (\cite{moses_instead_2020, moses_reverse-mode_2021}). 
This framework presents a few benefits. 
First, LLVM compilers allow us to use various programming languages in conjunction with this tool, opening up opportunities such as Julia and C++. 
Next, compile time optimization creates fast and reliable tools which are essential to generate code which can run at large scales on HPC.

For the work in section~\ref{sec:FiniteElements}, we use the MFEM library (\cite{anderson_mfem_2021}) due to its excellent scaling to HCP as well as various features such as partial assembly, sum factorization, native support for arbitrarily high order elements. 
Additionally, it is written in C++, which allows us to use Clang compilers (\href{https://llvm.org/}{https://llvm.org/}) to link with the Enzyme library for automatic differentiation.

\subsection{Automatic Differentiation at Different Scales (Defining Custom Vector-Jacobian Products)}
As mentioned previously, naive implementations of automatic differentiation can lead to inefficient solution processes or even be infeasible due to poor scaling and differences in hardware architectures. 
Recall the time stepping update given in Equation~(\ref{eq:discretedynamics}). 
The solution method is broken up into two main components, 1) the dynamical equation and 2) the form of time step update. 
In order to improve efficiency, we can break up the automatic differentiation into the same two distinct steps. 
We demonstrate this using a simple example. 

Consider Forward Euler time stepping. The update equation is given by
\begin{equation*}
	y_{k+1} = y_k + \Delta t_k f(y_k)\, .	
\end{equation*}
Taking the derivative with respect to $y_k$ and taking the inner product with a vector $\lambda$, we have
\begin{align*}
	\lambda \cdot \dfrac{\partial y_{k+1}}{\partial y_k} &= \lambda + \Delta t_k \lambda \cdot \dfrac{\partial f}{\partial y_k}
\end{align*}
Note that the derivative of the time stepping update depends on the derivative of the dynamics through the quantity $\lambda \cdot \dfrac{\partial f}{\partial y}$, but it does not depend on the exact form of this function directly. 
As a result, we can formulate modules for derivatives of time stepping operators with a separate function call to calculations which return $\lambda \cdot \dfrac{\partial f}{\partial y}$. 
In other words, changing the functional form of the physics does not change the structure of the derivative of the time stepping update, so that structure can be taken advantage of to improve efficiency. 
New time stepping formulations can also be implemented rapidly by using AD tools to take derivatives. 
There are additional reasons as to why this separation is beneficial when designing code for different hardware (HPC/GPU) that will become clear when discussing the specific case of Lagrangian hydrodynamics in Section~(\ref{sec:LagrangianHydrodynamics}).

Depending on the problem, there may be other levels at which scale separation can benefit code efficiency and modularity. 
This effect will be seen when discussing finite elements in Section~(\ref{sec:FiniteElements}).

\subsection{Static Quantities of Interest}
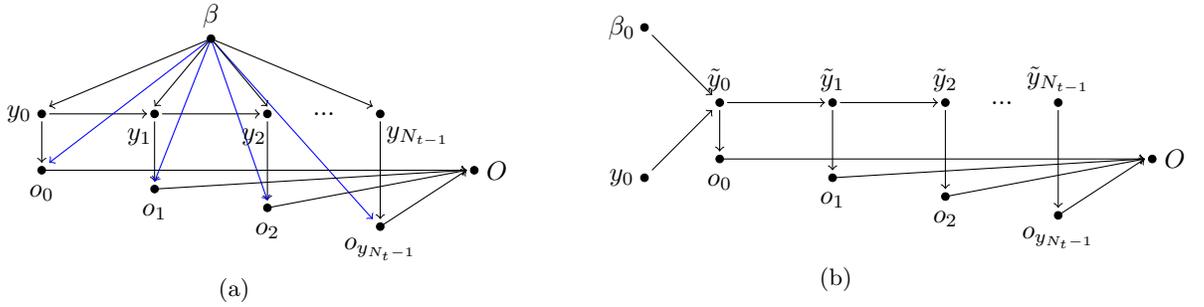
\begin{figure}[!hbt]
	\centering
	\begin{tikzpicture}
		\draw(0, 0) node[inner sep=0]{\begin{subfigure}[b]{.4\textwidth}
			\centering
			\begin{tikzpicture}
				\coordinate (AA) at (2.25,1);
				\coordinate (A) at (0,0);
				\coordinate (B) at (1.5,0);
				\coordinate (C) at (3.0,0);
				\coordinate (D) at (4.5,0);
				\coordinate (o0) at (0,-.75);
				\coordinate (o1) at (1.5,-1.);
				\coordinate (o2) at (3.0,-1.25);
				\coordinate (oyn) at (4.5,-1.5);
				\coordinate (O) at (5.75,-0.75);
				% \coordinate (E) at (6.0,0);
				\draw[->] ($(A) + (0.1, 0)$) -- ($(B) - (0.1, 0)$);
				\draw[->] ($(B) + (0.1, 0)$) -- ($(C) - (0.1, 0)$);
				\draw[->] ($(A) + (0.0, -0.1)$) -- ($(o0) + (0.0, 0.1)$);
				\draw[->] ($(B) + (0.0, -0.1)$) -- ($(o1) + (0.0, 0.1)$);
				\draw[->] ($(C) + (0.0, -0.1)$) -- ($(o2) + (0.0, 0.1)$);
				\draw[->] ($(D) + (0.0, -0.1)$) -- ($(oyn) + (0.0, 0.1)$);
				\draw[->,blue] (AA) -- ($(o0) + (0.10, 0.1)$);
				\draw[->,blue] (AA) -- ($(o1) + (0.0, 0.10)$);
				\draw[->,blue] (AA) -- ($(o2) + (0.0, 0.10)$);
				\draw[->,blue] (AA) -- ($(oyn) + (-0.10, 0.1)$);

				% \draw[->] ($(D) + (0.1, 0)$) -- ($(E) - (0.1, 0)$);

				\draw[->] ($(AA)$) -- ($(A) + (0.1, 0.1)$);
				\draw[->] ($(AA)$) -- ($(B) + (0.0, 0.10)$);
				\draw[->] ($(AA)$) -- ($(C) + (0.0, 0.10)$);
				\draw[->] ($(AA)$) -- ($(D) + (-0.1, 0.10)$);
				% \draw[->] ($(AA)$) -- ($(E) + (-0.1, 0.10)$);
				\draw ($(AA) + (0.0, 0.3)$) node {$\beta$};
				\draw [black,fill=black](AA) circle(0.05); 
				\draw ($(A) - (0.3, 0.0)$) node {$y_0$};
				\draw [black,fill=black](A) circle(0.05); 
				\draw ($(B) - (0.2, 0.3)$) node {$y_1$};
				\draw [black,fill=black](B) circle(0.05); 
				\draw ($(C) - (0.18, 0.3)$) node {$y_2$};
				\draw [black,fill=black](C) circle(0.05); 

				\draw ($(D) - (-0.5, 0.3)$) node {$y_{N_{t-1}}$};
				\draw [black,fill=black](D) circle(0.05);
				\draw ($(o0) - (0, 0.3)$) node {$o_0$};
				\draw [black,fill=black](o0) circle(0.05); 
				\draw ($(o1) - (0, 0.3)$) node {$o_1$};
				\draw [black,fill=black](o1) circle(0.05); 
				\draw ($(o2) - (0, 0.3)$) node {$o_2$};
				\draw [black,fill=black](o2) circle(0.05); 
				\draw ($(oyn) - (0, 0.3)$) node {$o_{y_{N_t - 1}}$};
				\draw [black,fill=black](oyn) circle(0.05); 
				% \draw ($(E) - (0, 0.3)$) node {$G$};
				% \draw [black,fill=black](E) circle(0.05); 
				\draw (3.75, 0) node {...};

				\draw ($(O) + (0.3, 0)$) node {$O$};
				\draw [black, fill=black](O) circle(0.05);
				\draw[->] ($(o0)$) -- ($(O) + (-0.1, 0.0)$);
				\draw[->] ($(o1)$) -- ($(O) + (-0.1, 0.0)$);
				\draw[->] ($(o2)$) -- ($(O) + (-0.1, 0.0)$);
				\draw[->] ($(oyn)$) -- ($(O) + (-0.1, 0.0)$);
			\end{tikzpicture}
			\caption{}\label{fig:static_quantities_a}
		\end{subfigure}
		};	
		\draw(8, 0) node[inner sep=0]{\begin{subfigure}[b]{.4\textwidth}
			\centering
			\begin{tikzpicture}
				\coordinate (AA) at (2.25,1);
				\coordinate (y0) at (-1, -1);
				\coordinate (beta) at (-1, 1);
				\coordinate (A) at (0,0);
				\coordinate (B) at (1.5,0);
				\coordinate (C) at (3.0,0);
				\coordinate (D) at (4.5,0);
				\coordinate (o0) at (0,-.75);
				\coordinate (o1) at (1.5,-1.);
				\coordinate (o2) at (3.0,-1.25);
				\coordinate (oyn) at (4.5,-1.5);
				\coordinate (O) at (5.75,-0.75);
				% \coordinate (E) at (6.0,0);
				\draw[->] ($(A) + (0.1, 0)$) -- ($(B) - (0.1, 0)$);
				\draw[->] ($(B) + (0.1, 0)$) -- ($(C) - (0.1, 0)$);
				\draw[->] ($(A) + (0.0, -0.1)$) -- ($(o0) + (0.0, 0.1)$);
				\draw[->] ($(B) + (0.0, -0.1)$) -- ($(o1) + (0.0, 0.1)$);
				\draw[->] ($(C) + (0.0, -0.1)$) -- ($(o2) + (0.0, 0.1)$);
				\draw[->] ($(D) + (0.0, -0.1)$) -- ($(oyn) + (0.0, 0.1)$);
				\draw[->] ($(y0) + (0.10, 0.1)$) -- ($(A) + (-0.1, -0.1)$);
				\draw[->] ($(beta) + (0.10, -0.1)$) -- ($(A) + (-0.10, 0.1)$);
				% \draw[->,blue] (AA) -- ($(o0) + (0.10, 0.1)$);
				% \draw[->,blue] (AA) -- ($(o1) + (0.0, 0.10)$);
				% \draw[->,blue] (AA) -- ($(o2) + (0.0, 0.10)$);
				% \draw[->,blue] (AA) -- ($(oyn) + (-0.10, 0.1)$);

				% \draw[->] ($(D) + (0.1, 0)$) -- ($(E) - (0.1, 0)$);

				% \draw[->] ($(AA)$) -- ($(A) + (0.1, 0.1)$);
				% \draw[->] ($(AA)$) -- ($(B) + (0.0, 0.10)$);
				% \draw[->] ($(AA)$) -- ($(C) + (0.0, 0.10)$);
				% \draw[->] ($(AA)$) -- ($(D) + (-0.1, 0.10)$);
				% \draw[->] ($(AA)$) -- ($(E) + (-0.1, 0.10)$);
				% \draw ($(AA) + (0.0, 0.3)$) node {$\beta$};
				% \draw [black,fill=black](AA) circle(0.05); 
				\draw ($(y0) - (0.3, 0.0)$) node {$y_0$};
				\draw [black,fill=black](y0) circle(0.05); 
				\draw ($(beta) - (0.3, 0.0)$) node {$\beta_0$};
				\draw [black,fill=black](beta) circle(0.05); 
				\draw ($(A) + (0.0, 0.3)$) node {$\tilde{y}_0$};
				\draw [black,fill=black](A) circle(0.05); 
				\draw ($(B) + (0.0, 0.3)$) node {$\tilde{y}_1$};
				\draw [black,fill=black](B) circle(0.05); 
				\draw ($(C) + (0.0, 0.3)$) node {$\tilde{y}_2$};
				\draw [black,fill=black](C) circle(0.05); 

				\draw ($(D) + (0.0, 0.3)$) node {$\tilde{y}_{N_{t-1}}$};
				\draw [black,fill=black](D) circle(0.05);
				\draw ($(o0) - (0, 0.3)$) node {$o_0$};
				\draw [black,fill=black](o0) circle(0.05); 
				\draw ($(o1) - (0, 0.3)$) node {$o_1$};
				\draw [black,fill=black](o1) circle(0.05); 
				\draw ($(o2) - (0, 0.3)$) node {$o_2$};
				\draw [black,fill=black](o2) circle(0.05); 
				\draw ($(oyn) - (0, 0.3)$) node {$o_{y_{N_t - 1}}$};
				\draw [black,fill=black](oyn) circle(0.05); 
				% \draw ($(E) - (0, 0.3)$) node {$G$};
				% \draw [black,fill=black](E) circle(0.05); 
				\draw (3.75, 0) node {...};

				\draw ($(O) + (0.3, 0)$) node {$O$};
				\draw [black, fill=black](O) circle(0.05);
				\draw[->] ($(o0)$) -- ($(O) + (-0.1, 0.0)$);
				\draw[->] ($(o1)$) -- ($(O) + (-0.1, 0.0)$);
				\draw[->] ($(o2)$) -- ($(O) + (-0.1, 0.0)$);
				\draw[->] ($(oyn)$) -- ($(O) + (-0.1, 0.0)$);
			\end{tikzpicture}
			\caption{}\label{fig:static_quantities_b}
		\end{subfigure}
		};	
	\end{tikzpicture}
	\caption{Graph diagrams for dynamical processes with static quantities of interest when (a) the variable is treated as global and affects every node in the graph and (b) the variable is treated pseudo-dynamically. }
\end{figure}

Not every variable of interest in topology optimization is the initial state $y_0$. In many cases, we want to study the sensitivity of our QOI with respect to some parameters which are global in scope. 
The graph structure of this system can be seen in Figure~\ref{fig:static_quantities_a}.

This optimization problem can be written as 

\begin{align*}
	O[y_0, \beta] &= \sum_{k = 0}^{N_t - 1} o(y_k, \beta, t_k) \Delta t_k
\end{align*}

Taking the derivative of this expression, can lead to complicated cross terms as every state of the solution $y_k$ is dependent on the parameter $\beta$.
One method to simplify the calculations is to instead treat $\beta$ as a time dependent parameter with a time derivative of $\dot{\beta} = 0$. 
As a result, we have the modified system of equations

\begin{align*}
	O[y_0, \beta_0] &= \sum_{k = 0}^{N_t - 1} o(y_k, \beta_k, t_k) \Delta t_k\, ,\\	
	\dot{y} &= f(y, \beta, t)\, , \\ 
	\dot{\beta} &= 0\, .
\end{align*}

While this change seems cursory at first, it allows us to define the adjoints in a much more straight-forward way. 
For example, we have
\begin{align*}
	\begin{bmatrix}
		\xi_y \\ 
		\xi_\beta
	\end{bmatrix} \cdot \begin{bmatrix}
		\dfrac{\partial y_{k+1}}{\partial y_k} && \dfrac{\partial y_{k+1}}{\partial \beta_k} \\ 
		\dfrac{\partial \beta_{k+1}}{\partial y_k} && \dfrac{\partial \beta_{k+1}}{\partial \beta_k}
	\end{bmatrix} &= \begin{bmatrix}
		\xi_y \\ \xi_\beta
	\end{bmatrix} \cdot \begin{bmatrix}
		\dfrac{\partial y_{k+1}}{\partial y_k}	&& \dfrac{\partial y_{k+1}}{\partial \beta_k} \\ 
		0 && 1
	\end{bmatrix} \\ 
	&= \begin{bmatrix}
		\xi_y^T \cdot \dfrac{\partial y_{k+1}}{\partial y_k} \\ 
		\xi_y^T \cdot \dfrac{\partial y_{k+1}}{\partial \beta_k} + \xi_\beta
	\end{bmatrix}
\end{align*}

From this, we see that we can utilize the same framework from previous sections. 
The only addition, other than keeping track of adjoints with respect to $\beta$, is the calculation of $\dfrac{\partial y_{k+1}}{\partial \beta_k}$. 
By accumulating the gradient through the time series, we now have the additional gradient information associated with $\dfrac{\partial O}{\partial \beta_0}$.

\subsection{Data Storage/Checkpointing}\label{sec:checkpointing}
When calculating adjoints for dynamic problems it is necessary to access the state data of every timestep from the forward solve during the reverse pass. 
This presents an issue in complex problems with many time steps and fine discretizations, as memory (RAM) can quickly become a limitation. 
For example, a typical hydrodynamic problem can involve solving for a system with a few million parameters. 
This system will may require $O(1M)$ time steps to fully resolve. Assuming we are storing each parameter in double precision, 

\begin{align*}
  \left(\frac{8 \text{ bytes}}{\text{parameter}} \right) \left(10^6 \text{ parameters} \right) \left( 10^6 \text{ time steps} \right) \approx 8 \text{ terabytes}	\, .
\end{align*}
This, unfortunately, is already exceeding the memory resources of current compute systems, and the problem becomes worse with more complex multiphysics and long duration simulations.

\begin{figure}[!hbt]
	\centering	
	\includegraphics[width=.5\textwidth]{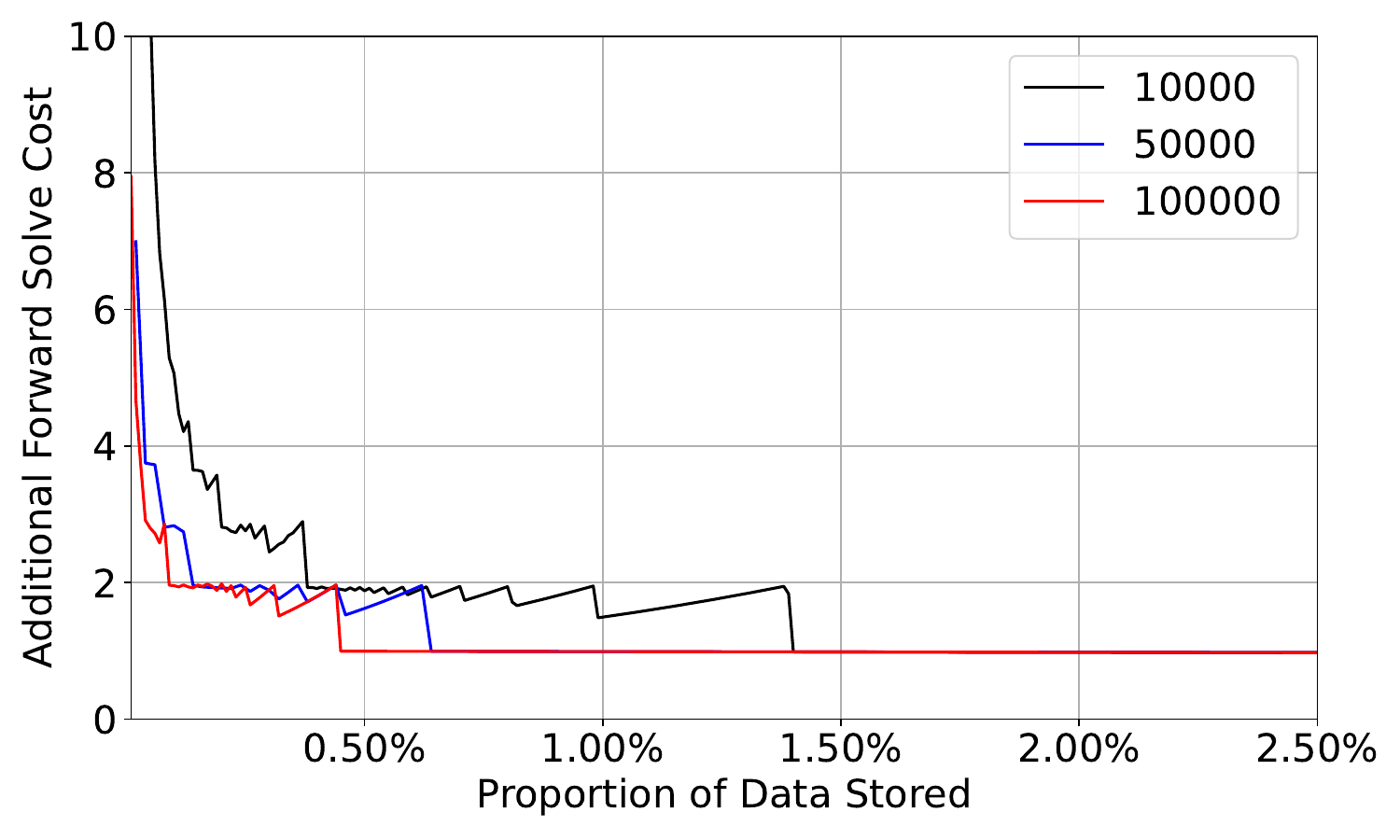}
	\caption{Amount of additional compute time added (measured in terms of the time of the forward solve) added given a memory budget written as a percentage of the total number of states from the forward solve. The different lines are for different numbers of total states from the forward solve.}\label{fig:checkpointing}
\end{figure}

One option for addressing this issue is to selectively store data to disk and load it into RAM when needed. 
This process is often inefficient in practice due to the relatively low bandwidth and high latency of disk reads and writes, leading to bad performance of both the forward and reverse solves.
Another option is the use a checkpointing scheme.  
The strategy here is to only save a relatively small number of simulation states in (relatively higher bandwidth) RAM.
On the reverse pass, the data for each state going back in time is required for the adjoint calculation.
When a state is reached that is not currently stored in RAM, the most recently saved checkpoint is fetched from RAM and the subsequent data states are recomputed by integrating forward in time again until the desired state is available.
Algorithms for constructing an optimal schedule for checkpointing and fetching to minimize the required number of recompute steps have been developed in:~\cite{griewank2000checkpoint}, which provides an optimal method when the number of time steps is known up front;~\cite{wang2009checkpoint}, which has good performance even when the number of steps is initially unknown (e.g., when the stable timestep evolves as the simulation progresses); and~\cite{hermann2020checkpoint} (and other references within), where different levels of available memory with varying costs and resources are considered.
Here we use the dynamic checkpointing algorithm from~\cite{wang2009checkpoint} to eventually allow for the possibility of varying timesteps, and we are currently only checkpointing to a single level of memory, namely, RAM.
The tradeoff between memory and compute time in our framework can be seen in Figure~\ref{fig:checkpointing}. 
In all cases shown, by storing only $1.5\%$ of the total data in an optimal manner, we incur an additional cost of about $1$ forward solve for recomputing the required states.  
The additional cost decays roughly linearly to zero in the limit where $100\%$ of the data is stored.

\subsection{What to do with gradients}
The above methods allow us to calculate the derivative of some objective function with respect to the initial (generalized) state. 
Access to this gradient, coupled with a gradient descent algorithm, allows us to minimize (or maximize) the given objective function.

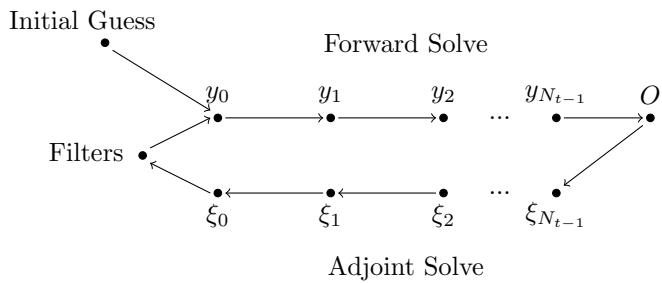
\begin{figure}[!hbt]
	\centering
	\begin{tikzpicture}
		\draw(0, 0) node[inner sep=0]{\begin{subfigure}[b]{\textwidth}
			\centering
			\begin{tikzpicture}
				\coordinate(IG) at (-1.5, 1);
				\coordinate (A) at (0,0);
				\coordinate (B) at (1.5,0);
				\coordinate (C) at (3.0,0);
				\coordinate (D) at (4.5,0);
				\coordinate (O) at (5.75,0.0);
				\coordinate (xiA) at (0,-1);
				\coordinate (xiB) at (1.5,-1);
				\coordinate (xiC) at (3.0,-1);
				\coordinate (xiD) at (4.5,-1);
				\coordinate (F) at (-1, -0.5);
				% \coordinate (E) at (6.0,0);
				\draw[->] ($(A) + (0.1, 0)$) -- ($(B) - (0.1, 0)$);
				\draw[->] ($(B) + (0.1, 0)$) -- ($(C) - (0.1, 0)$);
				\draw[->] ($(D) + (0.1, 0)$) -- ($(O) - (0.1, 0)$);
				\draw[->] ($(IG) + (0.1, -0.1)$) -- ($(A) + (-0.1, 0.1)$);

				\draw ($(IG) + (-0.3, 0.3)$) node {Initial Guess};
				\draw [black, fill=black](IG) circle(0.05);

				\draw ($(A) + (0.0, 0.3)$) node {$y_0$};
				\draw [black,fill=black](A) circle(0.05); 
				\draw ($(B) + (0.0, 0.3)$) node {$y_1$};
				\draw [black,fill=black](B) circle(0.05); 
				\draw ($(C) + (0.0, 0.3)$) node {$y_2$};
				\draw [black,fill=black](C) circle(0.05); 
				\draw ($(D) + (0.0, 0.3)$) node {$y_{N_{t-1}}$};
				\draw [black,fill=black](D) circle(0.05);
				\draw ($(F) + (-0.75, 0.05)$) node {Filters};
				\draw [black,fill=black](F) circle(0.05);
				\draw (3.75, 0) node {...};

				\draw ($(xiA) - (0.0, 0.3)$) node {$\xi_0$};
				\draw [black,fill=black](xiA) circle(0.05); 
				\draw ($(xiB) - (0.0, 0.3)$) node {$\xi_1$};
				\draw [black,fill=black](xiB) circle(0.05); 
				\draw ($(xiC) - (0.0, 0.3)$) node {$\xi_2$};
				\draw [black,fill=black](xiC) circle(0.05); 
				\draw ($(xiD) - (0.0, 0.3)$) node {$\xi_{N_{t-1}}$};
				\draw [black,fill=black](xiD) circle(0.05);
				\draw (3.75, -1) node {...};

				\draw[->] ($(xiB) - (0.1, 0)$) -- ($(xiA) + (0.1, 0)$);
				\draw[->] ($(xiC) - (0.1, 0)$) -- ($(xiB) + (0.1, 0)$);
				\draw[->] ($(O) - (0.1, 0.1)$) -- ($(xiD) + (0.1, 0.1)$);
				\draw[->] ($(xiA) + (-0.1, 0.)$) -- ($(F) - (-0.1, 0.1)$);
				\draw[->] ($(F) + (0.1, 0.1)$) -- ($(A) - (0.1, 0.0)$);
				\draw ($(O) + (0.0, 0.3)$) node {$O$};
				\draw [black, fill=black](O) circle(0.05);

				\draw (2.5, 1) node {Forward Solve};
				\draw (2.5, -2) node {Adjoint Solve};
			\end{tikzpicture}
		\end{subfigure}
		};	
	\end{tikzpicture}
	\caption{Graph outline of the optimization cycle starting with an initial guess $y_0$. The solution is propagated from $y_0$ to $O$. If the objective is not sufficiently converged, then the bottom path is taken to calculate the gradient. Then, the filters can be applied then the initial guess can be updated. This process can be iterated until convergence.}\label{fig:optimization_loop}
\end{figure}

For example, given an initial state $y_0$, we can calculate the forward solve to evaluate $O(y_0)$, then the adjoint solve to evaluate $\dfrac{\partial O}{\partial y_0}$. We then update the initial state with
\begin{align*}
	y_0 \leftarrow y_0 - \alpha \dfrac{\partial O}{\partial y_0} \, .
\end{align*}
An outline of this problem structure can be seen in Figure~(\ref{fig:optimization_loop}).

%\subsection{Problem Structure for Optimization}
Many of these optimization problems follow a similar structure (i) problem setup which defines physics and initial conditions, (ii) forward pass which advances time-stepping algorithm, (iii) objective function which maps final state(s) to a scalar value which we wish to minimize, and (iv) adjoint calculation via reverse pass where the gradient is accumulated.

\subsection{Example: Optimizing Multi-Particle Systems}\label{sec:ParticleDynamics}
For an introductory example, we introduce a system of interacting particles with two-body interactions. We label particles with degrees of freedom
\begin{align*}
	\text{Position} &- \mathbf{x}\\	
	\text{Velocity} &- \mathbf{v} \\ 
	\text{Charge} &- q
\end{align*}
Each particle is being acted upon by Coulomb and gravitational forces of the form
\begin{align*}
	\mathbf{F}_i = 	\sum_{j \neq i} \frac{q_i q_j \mathbf{r}_{ij}}{\| \mathbf{r}_{ij}\|^3} - g \mathbf{e}_2\, ,
\end{align*}
where the subscripts $i$ and $j$ indicate individual particles, $\mathbf{r}_{ij} = \mathbf{x}_i - \mathbf{x}_j$, $g$ is the gravitational constant, and $\mathbf{e}_2$ is the direction of the gravitational force. 
The dynamics of the system are given by
\begin{align*}
	\dot{\mathbf{x}}_i &= \mathbf{v}_i \\ 
	\dot{\mathbf{v}}_i &= \mathbf{F}_i \\ 
	\dot{q}_i &= 0
\end{align*}

Particles are initialized at random (non-overlapping) locations with zero initial velocity and uniform charge $q_i = 1.0$. 
The forward trajectory can be found by composing the physics with a time integration scheme. We choose a 4th order Runge-Kutta scheme summarized by
\begin{align}
	\begin{split}\label{eq:rk4}
		y_{k+1} &= y_k + \frac{\Delta t}{6}(k_1 + 2 k_2 + 2 k_3 + k_4)\, , \\
		k_1 &= f(t_k, y_k)\, , \\ 
		k_2 &= f \left( t_k + \frac{\Delta t}{2}, y_k + \frac{\Delta t}{2} k_1 \right)\, , \\ 
		k_3 &= f \left(t_k + \frac{\Delta t}{2}, y_k + \frac{\Delta t}{2} k_2 \right) \, , \\
		k_4 &= f \left( t_k + \Delta t, y_k + \Delta t k_3 \right) \, .
	\end{split}
\end{align}

%\subsection{Objective Function}
As an illustrative objective function, let's define a quantity of interest of the form
\begin{align*}
	O = \frac{1}{2}\sum_{i = 0}^{N_p - 1} \| \mathbf{x}_i^f - \mathbf{x}_i^* \|^2	
\end{align*}
where $\mathbf{x}_i^f$ is the final position of the particle, $\mathbf{x}_i^* = R (\cos \left( \frac{2 \pi i}{N_p} \right) \mathbf{e}_1 - \sin \left(\frac{2 \pi i }{N_p} \right) \mathbf{e}_2)$ and $R$ is the radius of the circle we wish to target.
The goal is to minimize the above function with respect to the initial state variables. We will specifically consider the case where we are only allowed to control the initial velocity of each particle and not the position or charge.
The initial derivative of the objective function can be found either through manual calculation or by using automatic differentiation.

%\subsection{Adjoint Calculation}
The adjoint of this integration scheme can be summarized as 
\begin{align*}
	m &= \xi \cdot y_{k+1} \\ 
	\bar{m} &= 1 \\ 
	\bar{y}_{k+1} &= \bar{m} \xi \\ 
	\bar{k}_4 &= \frac{\Delta t}{6}\bar{y}_{k+1} \\ 
	\bar{k}_3 &= \Delta t \bar{k}_4 \cdot \dfrac{\partial f}{\partial y}(t_k + \Delta t, y_k + \Delta t k_3) + \frac{\Delta t}{3}\bar{y}_{k+1} \\ 
	\bar{k}_2 &= \frac{\Delta t}{2}\bar{k}_3 \cdot \dfrac{\partial f}{\partial y} \left( t_k + \frac{\Delta t}{2}, y_k + \frac{\Delta t}{2} k_2 \right) + \frac{\Delta t}{3}\bar{y}_{k+1} \\ 
	\bar{k}_1 &= \frac{\Delta t}{2}\bar{k}_2 \cdot \dfrac{\partial f}{\partial y} \left( t_k + \frac{\Delta t}{2}, y_k + \frac{\Delta t}{2} k_1\right) + \frac{\Delta t}{6} \bar{y}_{k+1} \\ 
	\bar{y}_k &= \bar{y}_{k+1} + \bar{k}_4 \cdot \dfrac{\partial f}{\partial y}\left(t_k + \Delta t, y_k + \Delta t k_3 \right) + \bar{k}_3 \cdot \dfrac{\partial f}{\partial y}\left( t_k + \frac{\Delta t}{2}, y_k + \frac{\Delta t}{2} k_2 \right) \\ &+ \bar{k}_2 \cdot \dfrac{\partial f}{\partial y}\left( t_k + \frac{\Delta t}{2}, y_k + \frac{\Delta t}{2} k_1 \right) + \bar{k}_1 \cdot \dfrac{\partial f}{\partial y}\left( t_k, y_k \right)\, .
\end{align*}

By calculating $\bar{y}_k$ using the above algorithm, we have the incremental adjoint. Note that the above can be implemented manually, or calculated using an automatic differentiation package, defining custom gradients for the adjoints of the physics. By accumulating the incremental adjoint from the final time step to the initial time step using methods described in~\ref{sec:graphnetworkapproach}, we obtain the gradient with respect to the objective function.
Using the methods described in the previous sections, we begin with a set of initial states
\begin{align*}
	\mathbf{x}_i\, ,\mathbf{v}_i\, , q_i
\end{align*}
We can simulate the forward problem by using the 4th Order Runge-Kutta described in equations~\ref{eq:rk4}. 
In the solve, we cache all of the states in the forward solve. We then calculate the objective function and its initial gradient (with respect to the final state). 
Then, apply the incremental adjoint to march backwards in time (using the cached data) until the initial time step. 
At the end of this process, we have the gradient of the objective function with respect to the initial positions, velocities, and charges. 
We can then use an optimization algorithm, in our case conjugate gradient descent, to update the initial conditions. 
Because we only wish to optimize with respect to the initial velocities, we only update those quantities and zero out the perturbations of the positions and charges. 
This process can be repeated until convergence.

\begin{figure}[!hbt]
	\centering
	\begin{tikzpicture}
		\draw(0, 0) node[inner sep=0]{\begin{subfigure}[b]{.4\textwidth}
			\centering
			\includegraphics[width=\textwidth]{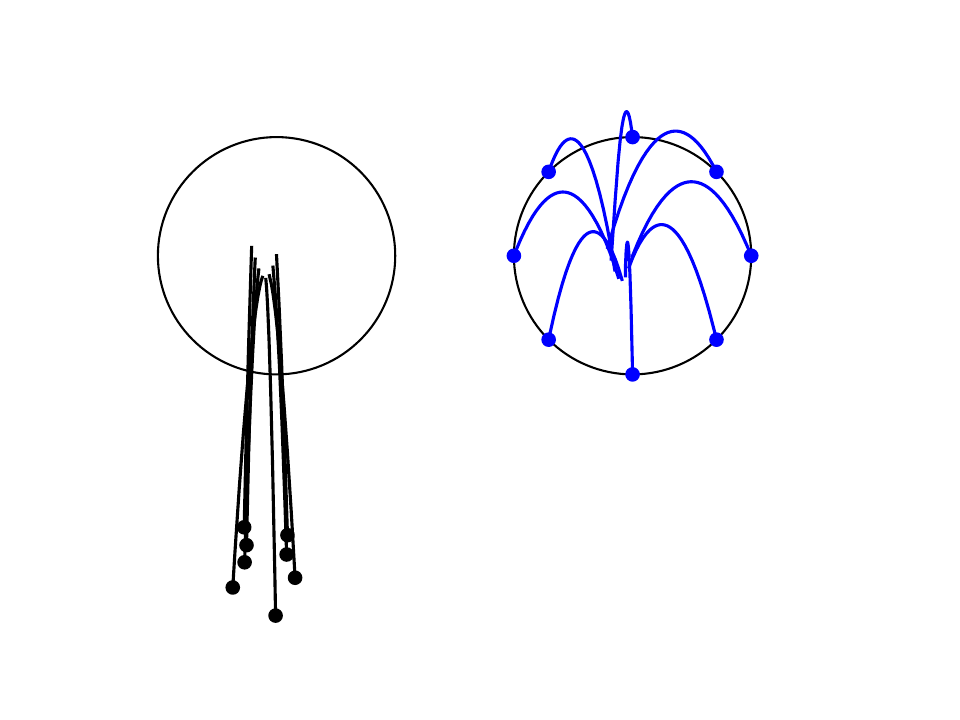}
			\caption{}
			\label{fig:particle_configs}
		\end{subfigure}
		};	
		\draw(6, 0) node[inner sep=0]{\begin{subfigure}[b]{.4\textwidth}
			\centering
			\includegraphics[width=\textwidth]{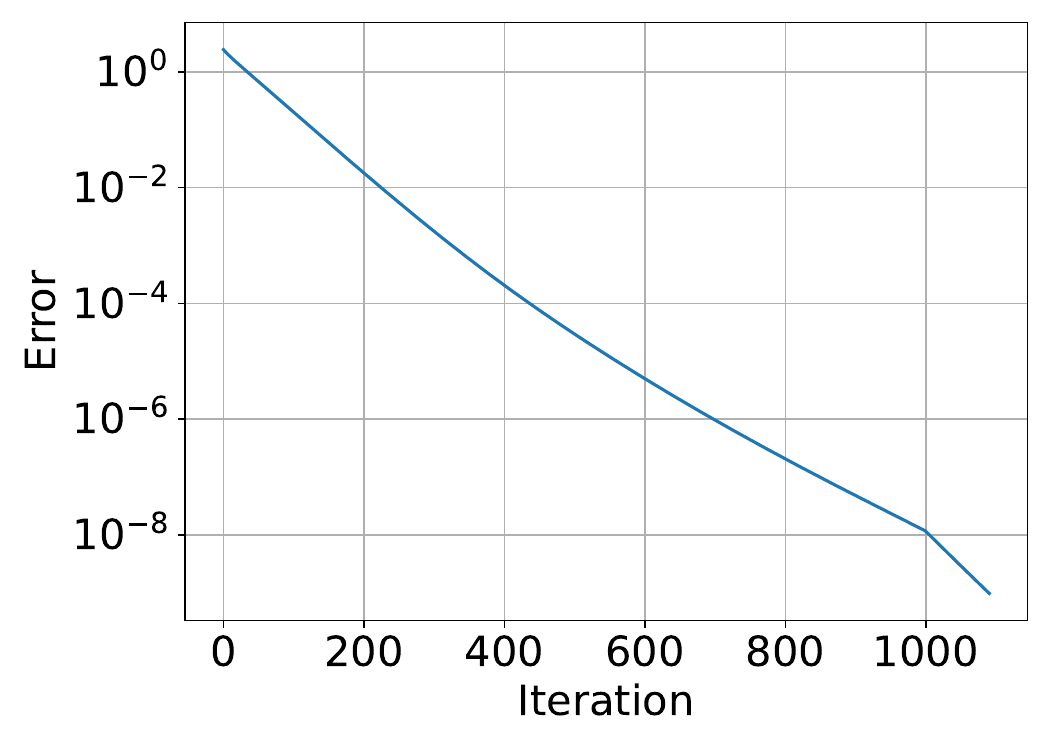}
			\caption{}
			\label{fig:particle_error}
		\end{subfigure}
		};	
		% Additional Writings
		\draw(-1.4, 2.1) node[inner sep=0]{\footnotesize Original};
		\draw(0.85, 2.1) node[inner sep=0]{\footnotesize Optimized};
	\end{tikzpicture}
	\caption{(a) Plots of the unoptimized (black) and optimized solutions (blue). (b) Plot of the error in the conjugate gradient solve demonstrating linear slope in the log-y scale.}

\end{figure}

While optimization of the positions of a system of particles seems straightforward, it actually illuminates many features of dynamic optimization. 
Inspection of the numerical solution of this particular particle system gives chaotic results in the informal sense that seemingly small changes may yield significantly different outcomes. 
We have found that there is often an interesting sweet spot in such hyperbolic systems wherein usage of gradient information is exceptionally meaningful; systems which are either very chaotic or not chaotic at all are rather boring from a gradient optimization perspective, whereas systems which are somewhat chaotic can be readily controlled employing gradient optimization.
While we have no formal proof of this behavior, we suspect that this trend is actually quite general.

\section{Finite Element Method discretization of Lagrangian Hydrodynamics}\label{sec:FiniteElements}
In this section, we apply the techniques discussed in the previous section to develop a method for computing adjoints of the equations of Lagrangian Hydrodynamics. We focus on a particular high-order discretization method which is described in~\cite{dobrev_high-order_2012}.
This requires special considerations for features like artificial viscosity and time step estimates. 
We review the spatial discretization of these equations into finite element bases, and the computation of the adjoint of each of these terms. 
\subsection{Summary of Equations}
The equations of Lagrangian Hydrodynamics (see~\cite{dobrev_high-order_2012, harlow_fluid_1971} for considerably more background and detail) describe the flow of continuous matter under the action of extreme pressures and energy deposition. 
The differential forms of the equations of motion, as defined in a Lagrangian reference frame are:
\begin{align}
	\text{Momentum conservation}: \quad \quad \rho \dfrac{d v}{dt} &= \nabla \cdot \sigma \, , \\ 
	\label{eq:diff_mass_cons}\text{Mass conservation}: \quad \quad \frac{1}{\rho} \dfrac{d \rho}{dt} &= - \nabla \cdot v \, , \\ 
	\text{Energy conservation}: \quad \quad \rho \dfrac{d e}{d t} &= \sigma : \nabla v \, , \\ 
	\text{Equation of motion}: \quad \quad \dfrac{d x}{d t} &= v\, , \\ 
	\label{eq:eos}\text{Stress Relation}: \quad \quad \sigma &= - p I + q \, ,
\end{align}
where $v$ is the material velocity, $e$ is the internal energy per unit density, $\rho$ is the current density, $p$ is the pressure calculated through an equation of state (EOS), and $q$ is an artificial viscosity.
Additionally, $\nabla$ is the spatial differential operator with respect to the current configuration.
The artificial viscosity regularizes strong shocks which otherwise would have a thickness significantly below the spatial resolution of the computational mesh~\cite{lew_artificial-viscosity_2001}. 
This topic will be discussed in more detail in Section~\ref{sec:artificial_viscosity}.

In semi-discretized form, these equations take the form
\begin{align}
	\text{Momentum conservation}: \quad \quad \mathbf{M}_v \dfrac{d \mathbf{v}}{dt} &= - \mathbf{F}\cdot \mathbf{1}\, , \\ 
	\text{Energy conservation}: \quad \quad \mathbf{M}_e \dfrac{d e}{dt} &= \mathbf{F}^T \cdot \mathbf{v} \, , \\	
	\text{Equation of motion}: \quad \quad \dfrac{d \mathbf{x}}{ dt} &= \mathbf{v}\, ,
\end{align}
where 
\begin{align}
	\mathbf{F}_{ij} = \int_{\Omega(t)} (\sigma : \nabla \phi_i^v) \phi_j^e dx\, ,	
\end{align}
$\mathbf{1}$ is a vector representing the constant 1, $\mathbf{M}_v$ is the mass matrix defined on the kinematic finite element space, and $\mathbf{M}_e$ is the mass matrix defined on the thermodynamic finite element space. 
Note that mass is conserved through the relation $\det(\nabla_\xi x) \rho= \det(\nabla_\xi X) \rho_0$, where $\xi$ are coordinates in a canonical frame, and $\rho$, $\rho_0$ are the densities in the current and reference configurations, respectively.
This enforces strong mass conservation point-wise and eliminates the need to explicitly solve for Equation~(\ref{eq:diff_mass_cons}). 
A consequence of this relation is that the mass matrices ($\mathbf{M}_v$ and $\mathbf{M}_e$) are independent of time.
We use a Mie-Guneisen equation of state of the form 
\begin{align*}
	p = \frac{\rho_0 C_0^2 \chi }{(1 - s \chi)^2}\left( 1 - \frac{\Gamma_0}{2} \chi \right) + \rho_0 \Gamma_0 e\, , \quad \quad \chi = 1 - \frac{\rho_0}{\rho}
\end{align*}
where $\Gamma_0$ is the Gruneisen parameter and $s$ is the linear Hugoniot slope coefficient.
\subsection{Artificial Viscosity}\label{sec:artificial_viscosity}
As stated in the previous section, artificial viscosity is critically important when considering shock propagation in order to prevent spurious oscillations without damping out the features of the shock. 
Various works, see~\cite{campbell_tensor_2001,lew_artificial-viscosity_2001,wilkins_use_1980}, discuss the formulations and derivations.
We use the particular form of artificial viscosity~\cite{wilkins_use_1980,dobrev_high-order_2012} given by 
\begin{align}
	\sigma_v = 0.75 \rho (\gamma_1 l c + \gamma_2 l |\Delta v|) H\left(\Delta v \right)	\text{sym}(\nabla v)
\end{align}
where $\gamma_i$ are strength parameters, $l$ is a length scale associated with an element ($l = l_0 \det(F)^{1/\text{dim}}$), $l_0$ is the initial length scale, $c$ is the wave speed in the element, $\nabla v = \text{tr} (\nabla v)l$ is the velocity jump across the element, and $H$ is the Heaviside function. 
This form has been has been demonstrated to sufficiently resolve the shock front while not overly damping the rest of the behavior; although we note that for truly high-order dissipation, which we will not consider here.
Additional limiting is required, such as the hyperviscosity treatment shown in~\cite{BELLOMALDONADO2020}. 
The choice of artificial viscosity, however, does not change the mathematical formulation, as the automatic differentiation tools allow for seamless transition between functional forms without having to recalculate derivatives.

One complication of this function in the context of calculating adjoints is the compression switch. 
In standard methods, we use a Heaviside function which is both non-differentiable, but also discontinuous. 
This is partially remedied by being multiplied by $\text{sym}(\nabla v)$, which raises this to being continuous but non-differentiable. 
We will discuss the issue of non-differentiability in a future work. 
For now, we remedy this by replacing all non-differentiability with a suitable smoothing function. 
In the case of the Heaviside function, we use a sigmoid scaled with the wave speed
\begin{align*}
	H(- x) \rightarrow \text{sigmoid}\left( -\frac{x}{h}\right) \, .
\end{align*}
Additionally, the absolute value in the quadratic term must replaced with a soft absolute value. This is done by
\begin{align*}
	|x| \rightarrow \text{softabs} (x, h) = \text{silu} \left( \frac{x}{h} \right) + \text{silu}\left( \frac{-x}{h} \right)
\end{align*}
where silu is the sigmoid linear unit and $h$ is a length scale.

In both these cases we choose $ h = 0.2 c$ to properly scale the transition in the smoothed regions to properly account for the behavior yet maintain differentiability where $c$ is a representative wave speed.

\subsection{Finite elements}
While the derivations for adjoint methods described in~\ref{sec:graphnetworkapproach} are appropriate for finite element methods, additional care must be taken when passing derivatives through the forcing functions. 
A standard dynamic finite element problem will include a forcing function, mass matrix, and time stepping algorithm.
% A standard discretization scheme for a dynamic finite element problem will include the following:
% \begin{itemize}
	% \item Forcing Function \\
		The forcing function often takes the form
		\begin{align} \label{eq:forcing}
			F_i = \int_\Omega f(x, v(x), \nabla v(x)) \cdot \phi_i(x) dx\, ,
		\end{align}
		where $\Omega$ is the computational domain, $\phi_i(x)$ are shape functions and $f$ is a local forcing function. 
		The multiplication in the final part of this integral is often done with either the shape functions $\phi_i$ or their gradients $\nabla \phi_i(x)$. 
By using the Galerkin formulation, our solution fields are expanded in terms of the shape functions as
\begin{align}
	v(x) = \sum_i v_i \phi_i(x)\, .
\end{align}

		The product in equation~\ref{eq:forcing} is linear with respect to the data $f$ and can be calculated in parallel. 
		Much work has been done in implementing this product in HPC (High Performance Computing) and GPUs (Graphics Processing Units).
	% \item Mass Matrix\\ 

		The mass matrix is often of the form
		\begin{align}
			M_{ij} = \int_\Omega \rho (x) \phi_i(x) \cdot \phi_j(x) dx\, ,
		\end{align}
		is a bilinear form where $\rho$ is the density field. 
		A complication is that this quantity is rarely utilized directly; instead we compute the action of its inverse ($M_{ij}^{-1} v_j$). 
		As a result, special care must be taken when differentiating.
	% \item Time stepping
% \end{itemize}

The final component is the time-stepping algorithm used to map the state at a particular timestep to the state at thee next timestep.
Recall from equation~\ref{eq:discretedynamics}, this operation is formed from the composition of the particular algorithm used (Euler, Runge-Kutta, etc.) and the dynamics (time derivative terms formulated using the above two calculations).
The integrator we choose largely depends on characteristics of the dynamics we wish to satisfy; for example, energy conservation, condition, and stability.

\subsection{Adjoints of Forcing}
To take the adjoint product of a force, we need to calculate
\begin{align*}
	\xi \cdot \dfrac{\partial \mathbf{F}}{\partial v_j} &= \sum_i \xi_i \dfrac{\partial F_i}{\partial v_j} \\ 
	&= \sum_i \xi_i \dfrac{\partial}{\partial v_j}\left( \int_\Omega f(x, v(x), \nabla v(x))\cdot \phi_i(x) dx \right) \\ 
	&= \sum_i \xi_i \int_\Omega \left( \dfrac{\partial f}{\partial v} \cdot  \dfrac{\partial v}{\partial v_j} + \dfrac{\partial f}{\partial \nabla v} \cdot \dfrac{\partial \nabla v}{\partial v_j} \right)\cdot \phi_i(x) dx \\ 
	&= \sum_i \xi_i \int_\Omega \left( \dfrac{\partial f}{\partial v} \cdot \phi_j  + \dfrac{\partial f}{\partial \nabla v} \cdot \nabla \phi_j \right)\cdot \phi_i(x) dx \\ 
	&= \int_\Omega \left( \dfrac{\partial f}{\partial v} \cdot \phi_j(x) + \dfrac{\partial f}{\partial \nabla v} \cdot \nabla \phi_j(x)\right) \cdot \xi(x) dx \, .
\end{align*}

By using the Galerkin expansion, we can project the finite element discretization of the discrete adjoint vector onto the shape functions and interpret that product as an adjoint field. 
Additionally, the calculation of the adjoint load takes a similar form to the calculation of the force load and, therefore, similar machinery can be used to calculate both.

\subsection{Adjoints of Mass Matrices}
The mass matrix plays an important role in a dynamical system. 
In some cases, we would like to calculate the sensitivity of the objective with respect to the initial density field. 
This idea is similar to the canonical problem in topology optimization of deciding where to place mass to minimize compliance. 
Because of this, it is important to consider the adjoint of the density field.

Let the density be a function of the state ($\rho = \rho(v(x))$). Then, when calculating the derivative, we have
\begin{align*}
	\sum_i \xi_i \dfrac{\partial}{\partial v_j}\left(M_{ik}^{-1}\right) &= - \sum_{ilm} \xi_i M_{il}^{-1} \dfrac{M_{lm}}{\partial v_j}M_{mk}^{-1}
\end{align*}
where we use that $M$ is invertible and $M^{-1} M = I$. Passing the derivative inside the integral, we have
\begin{align*}
	\dfrac{\partial}{\partial v_j} \int_\Omega \rho(v(x)) \phi_l(x) \cdot \phi_m(x) dx &= \int_\Omega \dfrac{\partial \rho}{\partial v} \phi_j(x) \phi_l(x) \cdot \phi_m(x) dx \, .
\end{align*}

Note that the above quantity is a third order tensor. 
While it will retain the same sparsity of the system as a whole, we would still like to avoid constructing such a quantity. 
Note that this will only appear with two contractions. 
Generally, we will see this as 
\begin{align*}
	\sum_{lm} \xi_l \eta_m \int_\Omega \dfrac{\partial \rho}{\partial v} \phi_j(x) \phi_l(x) \cdot \phi_m(x) dx &= \int_\Omega \dfrac{\partial \rho}{\partial v} \phi_j(x) \left( \sum_l \xi_l \phi_l(x) \right) \cdot \left( \sum_m \eta_m \phi_m(x) \right) dx \\ 
	&= \int_\Omega \dfrac{\partial \rho}{\partial v} \phi_j(x) \left( \xi(x) \cdot \eta(x) \right) dx
\end{align*}

By constructing the action of the tensor rather than the full tensor itself, we can drastically improve performance and avoid storing matrices~\cite{vargas2022matrix}. 
This technique is common in matrix free methods and is essential when studying extremely large systems.

\iffalse

\begin{figure}[H]
	\begin{center}
		\includegraphics[width=0.5\textwidth]{figures/fig1/plot.pdf}
	\end{center}
\end{figure}

\begin{figure}[H]
	\begin{center}
		\def\firstcircle{ (0.0, 0.0) circle (1.5)}
		\def\secondcircle{(2.0, 0.0) circle (1.5)}
		\def\thirdcircle{ (1.0,-1.5) circle (1.5)}
		\def\rectangle{ (-1.5,-3.0) rectangle (3.5,1.0) }
		\colorlet{circle edge}{black}
		\colorlet{circle area}{blue!30}

		\tikzset{filled/.style={fill=circle area, draw=circle edge, thick},
		outline/.style={draw=circle edge, thick}}

		\begin{tikzpicture}
			\begin{scope}
				\fill[filled]  \firstcircle;
				\fill[white]  \secondcircle;
				\fill[white] \thirdcircle;
			\end{scope}
			\draw[outline] \firstcircle  node[left]  {$A$};
			\draw[outline] \secondcircle node[right] {$B$};
			\draw[outline] \thirdcircle  node[below] {$C$};
			\node[anchor=south] at (current bounding box.north) 
			{$P = (A - B) - C$};
		\end{tikzpicture}
		%
		\begin{tikzpicture}
			\begin{scope}
				\fill[filled]  \firstcircle;
				\fill[white]  \secondcircle;
				\fill[white] \thirdcircle;
			\end{scope}
			\draw[outline] \firstcircle  node[left]  {$A$};
			\draw[outline] \secondcircle node[right] {$B$};
			\draw[outline] \thirdcircle  node[below] {$C$};
			\node[anchor=south] at (current bounding box.north) 
			{$Q = (A - C) - (B - C)$};
		\end{tikzpicture}
	\end{center}
\end{figure}
\fi
\subsection{Filters}
A common difficulty when calculating adjoints in mechanics is the existence of checkerboarding modes (rapid oscillations due to weak convergence of the gradient with respect to the mesh)~\cite{bendsoe_topology_2004,allaire_shape_2002,cherkaev_variational_2000}. 
Although this is generally discussed in the context of density based optimizations, rapid oscillations are pervasive when calculating derivatives with respect to finite element spaces. 
A common solution to this problem is to apply a filter of some sort to regularize the final adjoint before adjusting parameters for optimization. 
In many cases, addressing checkerboarding can also yield more practical and manufacturable designs for engineering requirements. 
In the following sections, we will use a modified mass matrix as a filter. 
This was chosen for its simplicity as many other filters could be considered.

Let be the derivative of the adjoint be $\dfrac{\partial o}{\partial y_0}$. To avoid propagating the checkerboarding pattern to the update step, we calculate 
\begin{align*}
	\Delta y = \mathbf{M}_\mathbf{1}^{-1} \dfrac{\partial o}{\partial y_0}	
\end{align*}
where $\mathbf{M}_{\mathbf{1},ij} = \int_\Omega \phi_i(x) \phi_j(x) dx$ is a mass matrix of unit density. The filtered update can be used in place of the original gradient for the purposes of optimization.

\section{Suppression of Richmyer-Meshkov Instability induced jetting}\label{sec:LagrangianHydrodynamics}
As noted in the introduction, RMI plays an important role in many scientific and engineering applications. Depending on the use case, either enhancement or suppression of the RMI jet may be desired. 
In this section, we develop a computational model of RMI, multi-objective functions to describe what we are after in terms of stable shock acceleration, the adjoint computation, some specific discussion of tracer particles practically needed to implement our specific objective function, and results showing RMI suppression.

\subsection{Forward Pass}
The domain is separated into two regions. The left is high density ($\rho_0 = 10$) and the right is low density ($\rho_0 = 1$). 
A subset of the left domain ($\Omega_1 \subset \Omega \text{  s.t. } X < 1$) is the controllable domain and is initialized to a higher internal energy state $e(X \in \Omega_1,0) = 0.15$ with the rest $e(X \notin \Omega_1, 0) = 0.0$. 
The top, left, and bottom boundaries allow for sliding boundary conditions while the right boundary is completely free.

A shock wave is generated due to the internal energy of the left side being higher than the rest of the domain. 
This causes a high pressure region with a sharp interface against a low pressure region, resulting in a force along the interface. 
As the high energy region expands, it generates a shock wave which propagates through the high density region. 
Once it hits the interface, baroclynic torque is generated due to the misalignment between the pressure gradient and the density gradient. 
This torque causes the system to evolve in such a way that the interface will invert itself and continue to grow. 
As before, we use the 4th order Runge-Kutta for time integration and solve until $t=7$.

The goal is to design the profile of the internal energy in $\Omega_1$ to minimize the RMI jet length at a particular time.

\subsection{Objective Function}
In order to minimize the jet length, we need a metric which takes the state of the system as an input and returns a single scalar value. 
By minimizing said functional, we obtain better results (i.e. reducing the jet length). 
Additionally, we would like to push the interface and not remove all accelerations (i.e. we want the optimizer to avoid the trivial solution of no acceleration). 
As a result, we would also like to include a term in the objective that increases the velocity of the interface while decreasing the objective value. In order to accomplish this, we introduce Lagrangian tracer particles. 
These are virtual particles which are linked to particular material points. In an HPC environment, there are complexities defining these in a method that allows for adjoint calculations. 
This is due to the question of ownership of data and how to communicate information about tracers across multiple threads. 
The exact implementation of these objects, however, is outside the scope of this manuscript.

\begin{figure}[hbt!]
	\centering
	\begin{tikzpicture}
		\draw(0, 0) node[inner sep=0, anchor=north west]{\begin{subfigure}[b]{.6\textwidth}
			\centering
			\includegraphics[width=\textwidth]{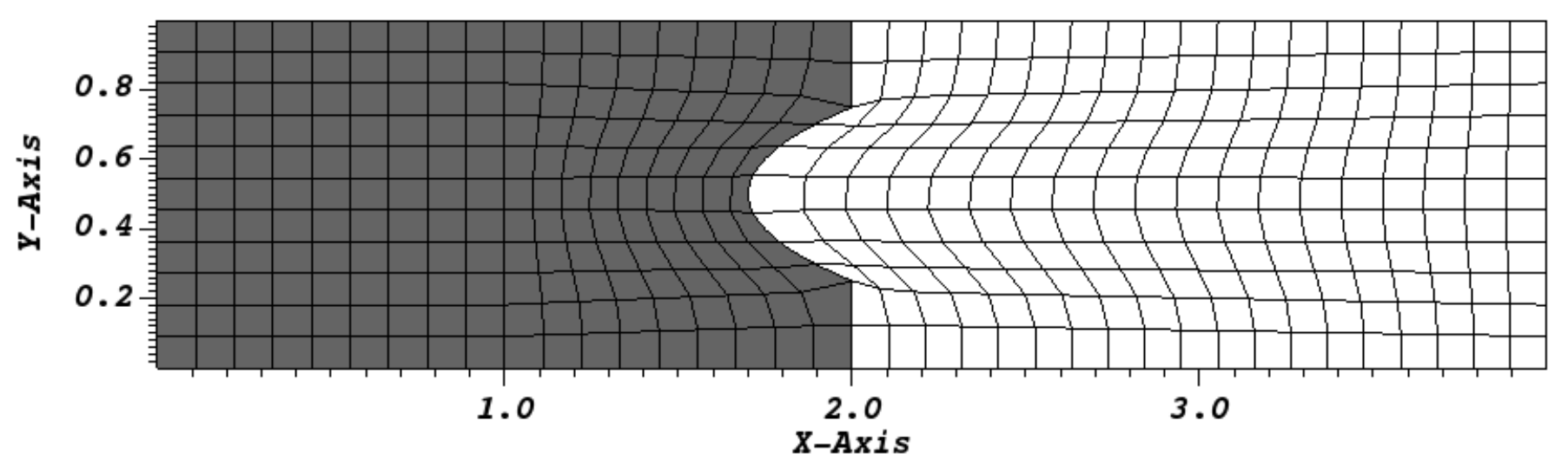}
		\end{subfigure}
		};	
		\draw(0, -3) node[inner sep=0, anchor=north west]{\begin{subfigure}[b]{.6\textwidth}
			\centering
			\includegraphics[width=\textwidth]{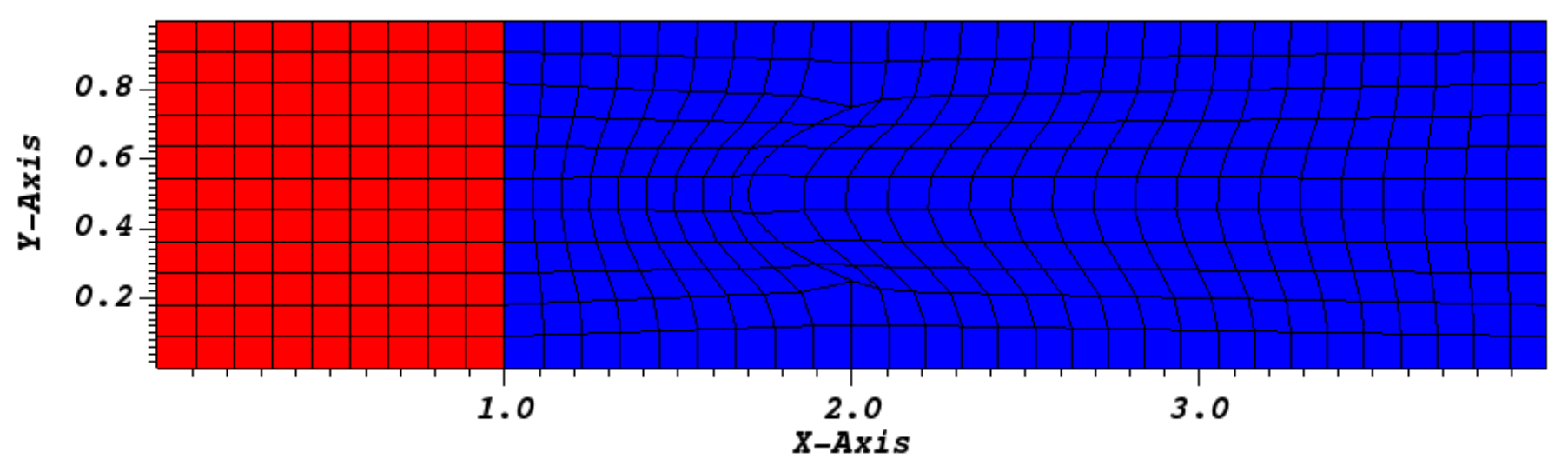}
			% \caption{}
		\end{subfigure}
		};	
		% Additional Writings
		\draw(-0.3, 0.3) node[anchor=north west, inner sep=0]{\footnotesize Density};
		\draw(-0.3, -2.7) node[anchor=north west, inner sep=0]{\footnotesize Energy};
		\draw [dashed, black!70] (2.92, 0.3) -- (2.92, -6.0);
		% \draw [dashed, blue!50] (5.12, -3.6) -- (1.42, -4.95);
		% \draw [dashed, red] (10.12, -3.6) -- (1.45, -5.1);
		\draw [decorate,decoration={brace,amplitude=5pt,raise=4ex}]
		(.92,-.5) -- (4.9,-.5) node[midway,yshift=3em]{\footnotesize High Density};
		\draw [decorate,decoration={brace,amplitude=5pt,raise=4ex}]
		(4.9,-.5) -- (8.9,-.5) node[midway,yshift=3em]{\footnotesize Low Density};
		\draw [decorate,decoration={brace,amplitude=5pt,mirror,raise=4ex}]
		(0.92,-4.9) -- (2.85,-4.9) node[midway,yshift=-3em]{\footnotesize $\Omega_1$};
		\draw [red,fill=red](4.82,-.25) circle(0.05); 
		\draw (5.05, -0.25) node[inner sep=0]{\footnotesize 2};
		\draw [red,fill=red](4.82,-2.02) circle(0.05); 
		\draw (5.05, -2.02) node[inner sep=0]{\footnotesize 3};
		\draw [red,fill=red](4.28,-1.125) circle(0.05); 
		\draw (4.50, -1.125) node[inner sep=0]{\footnotesize 1};
	\end{tikzpicture}
	\caption{Domain and initial energy configuration for the RMI case study. The nodes (1, 2 ,3) indicate where the tracer particles are placed (used in Equation~\ref{eq:rmi_objective}).}\label{fig:rmi_config}

\end{figure}
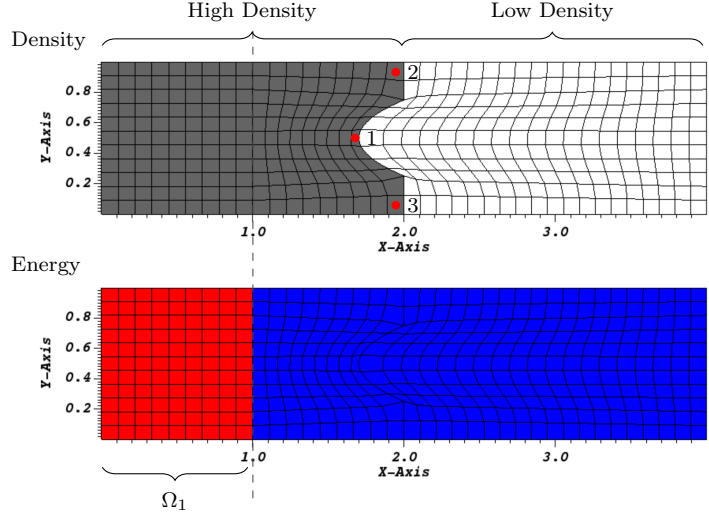

We introduce the objective
\begin{align}\label{eq:rmi_objective}
	O = \frac{1}{2} \lambda_1 (x_1 - x_\text{outer})^2 +  \frac{\lambda_2}{\delta + |v_\text{ave}|}	
\end{align}
where $x_i$ are the $x$ components of the deformation of particles $i$, $v_i$ are the $X$ components of the velocities of particles $i$, $x_\text{outer} = \text{ave}(x_2, x_3)$, $v_\text{ave} = \text{ave}(v_1, v_2, v_3)$, and $\lambda_i \geq 0$ are scaling factors. 
By inspection, it can be seen that this function is minimized when
\begin{align*}
	x_\text{outer} = x_1 \, , \quad \quad &\text{Flatten the interface} \\	
	v_\text{ave} = \pm \infty\, , \quad \quad &\text{Accelerate interface}
\end{align*}

when $\lambda_i > 0$.

In our particular case, we would like to push $v_\text{ave} \to +\infty$ as opposed to the other side. 
This is remedied by picking a suitable initial guess and using a local optimizer to push the solution towards one solution as opposed to the other. 
If using completely random initial conditions or a global minimizer, then a different objective function may be necessary.
Also, note that Lagrangian tracer particles use local data to construct their state, as a result, we can generally write that
\begin{align} \label{eq:tracer}
	\begin{split}
		\mathbf{x}_i &= \mathbf{x}_i(X, x, v, e) \, , \\ 
		\mathbf{v}_i &= \mathbf{v}_i(X, x, v, e) \, , \\ 
		e_i &= e_i(X, x, v, e) \, ,
	\end{split}
\end{align}
where the terms on the left are features of the tracer particles and the equations on the right are functions which project the global state vectors to the tracer data. 
\subsection{Adjoint Calculation}
As a result of the tracer definition in~(\ref{eq:tracer}), we can calculate the gradient of the objective as
\begin{align*}
	\dfrac{\partial O}{\partial X} &= \sum_i\left(  \dfrac{\partial O}{\partial \mathbf{x}_i} \cdot \dfrac{\partial \mathbf{x}_i}{\partial X} +  \dfrac{\partial O}{\partial \mathbf{v}_i }\cdot \dfrac{\partial \mathbf{v}_i}{\partial X} + \dfrac{\partial O}{\partial e_i} \dfrac{\partial e_i}{\partial X}\right) \, , \\ 
	\dfrac{\partial O}{\partial x} &= \sum_i\left(  \dfrac{\partial O}{\partial \mathbf{x}_i} \cdot \dfrac{\partial \mathbf{x}_i}{\partial x} +  \dfrac{\partial O}{\partial \mathbf{v}_i }\cdot \dfrac{\partial \mathbf{v}_i}{\partial x} + \dfrac{\partial O}{\partial e_i} \dfrac{\partial e_i}{\partial x}\right) \, , \\ 
	\dfrac{\partial O}{\partial v} &= \sum_i\left(  \dfrac{\partial O}{\partial \mathbf{x}_i} \cdot \dfrac{\partial \mathbf{x}_i}{\partial v} +  \dfrac{\partial O}{\partial \mathbf{v}_i }\cdot \dfrac{\partial \mathbf{v}_i}{\partial v} + \dfrac{\partial O}{\partial e_i} \dfrac{\partial e_i}{\partial v}\right) \, , \\ 
	\dfrac{\partial O}{\partial e} &= \sum_i\left(  \dfrac{\partial O}{\partial \mathbf{x}_i} \cdot \dfrac{\partial \mathbf{x}_i}{\partial e} +  \dfrac{\partial O}{\partial \mathbf{v}_i }\cdot \dfrac{\partial \mathbf{v}_i}{\partial e} + \dfrac{\partial O}{\partial e_i} \dfrac{\partial e_i}{\partial e}\right) \, . \\ 
\end{align*}
Note that the relations above tend to be simple, but we maintain all the terms for the sake of generality. 
Similar to previously discussed methods, we calculate the adjoints of the time integration scheme and of the physics separately and compose them in order to find the incremental adjoint.
By stepping backwards in time until the initial timestep, we can accumulate the adjoint of the objective with respect to the entire initial state $X, x, v, e$. 
From here, we mask the adjoint to only include the components in $\Omega_1$. 
\begin{figure}[!hbt]
	\centering
	\includegraphics[width=\textwidth]{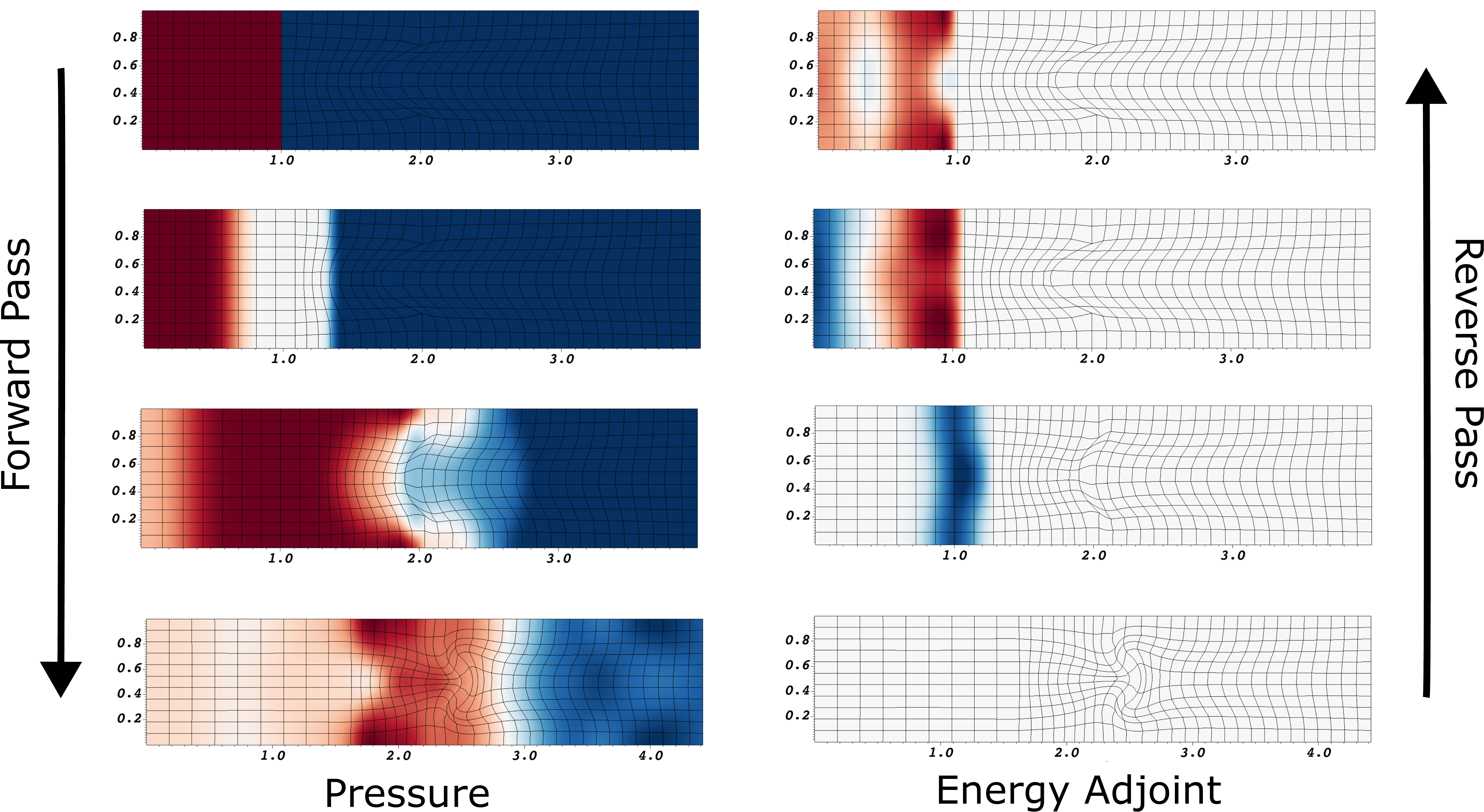}
	\caption{Visualization of pressure and adjoint fields through the solve. The top plots are solutions at time $t=0$ while the bottom are the final steps at $t = 7$. Following the result from the top left and working down, we integrate forward in time. Then from the bottom left to bottom right, we calculate the objective and initial gradient. Moving up, we integrate backwards in time to solve for the adjoint energy field to calculate the gradient. Note that the system has the density interface seen in Figure~\ref{fig:rmi_config}.}\label{fig:rmi_adjoint_loop}
\end{figure}

An example of the forward and adjoint loops are shown in Figure~\ref{fig:rmi_adjoint_loop}.
We plot the pressure evolution from the forward pass; however, we cache all the data needed to recreate the state using the checkpointing methods from Section~\ref{sec:checkpointing}.
Then, we calculate our objective function (see equation~\ref{eq:rmi_objective}) and its gradient with respect to the final state. 
Performing adjoint calculations, we step backwards in time, recreating the state as needed using the checkpoints.
We visualize the adjoint of the energy field (masked by the domain of control $\Omega_1$) as we step backward in time. 
At the final time step of the adjoint solve ($t=0$) we are left with the gradient of our objective with respect to the initial state within $\Omega_1$.
Visualization of the adjoint fields is extremely useful when trying to understand the system's sensitivities.
Additionally, these fields provide vital information which designers can use to improve designs ad hoc.

\subsection{Results}\label{sec:RMI_Optimization}
\begin{figure}[!hbt]
% \begin{figure}[H]
	\centering
	\begin{tikzpicture}
		\draw(0, 0) node[inner sep=0, anchor=north]{\begin{subfigure}[b]{.33\textwidth}
			\centering
			\includegraphics[width=\textwidth]{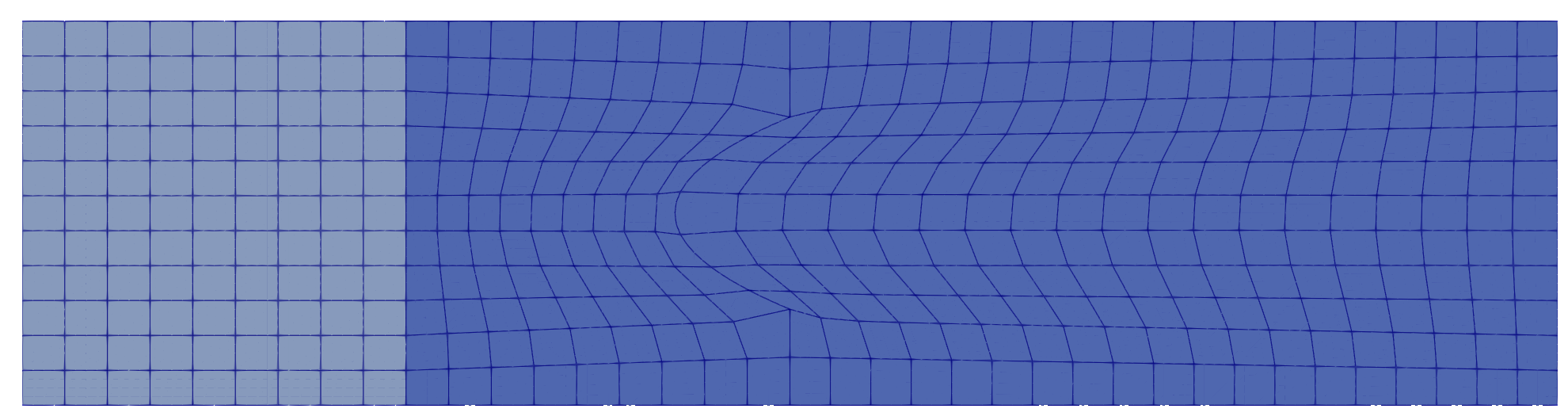}
		\end{subfigure}
		};	
		\draw(5, 0) node[inner sep=0, anchor=north]{\begin{subfigure}[b]{.33\textwidth}
			\centering
			\includegraphics[width=\textwidth]{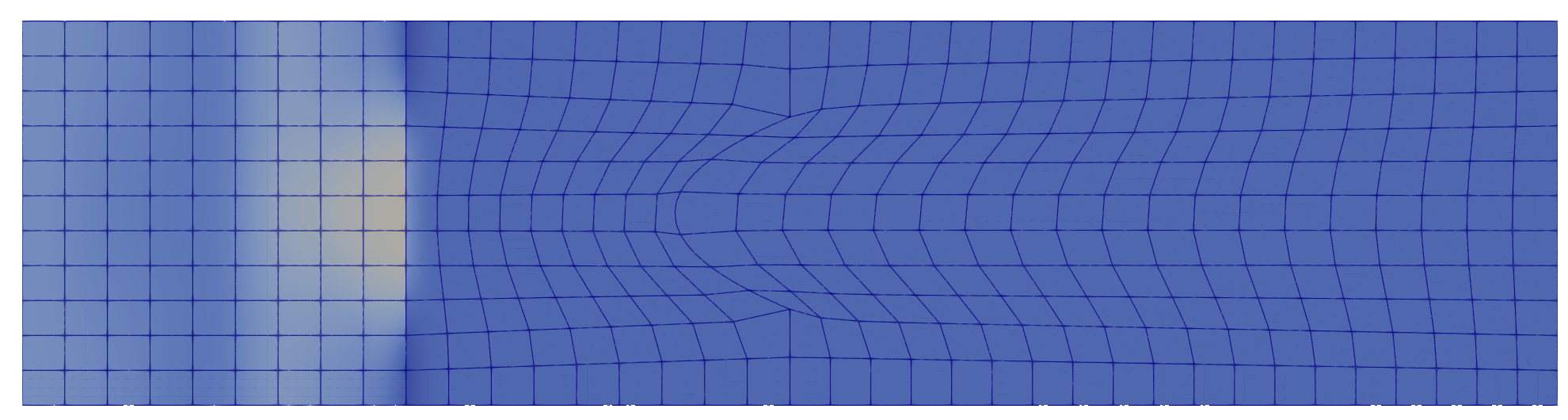}
			\caption{}
		\end{subfigure}
		};	
		\draw(10, 0) node[inner sep=0, anchor=north]{\begin{subfigure}[b]{.33\textwidth}
			\centering
			\includegraphics[width=\textwidth]{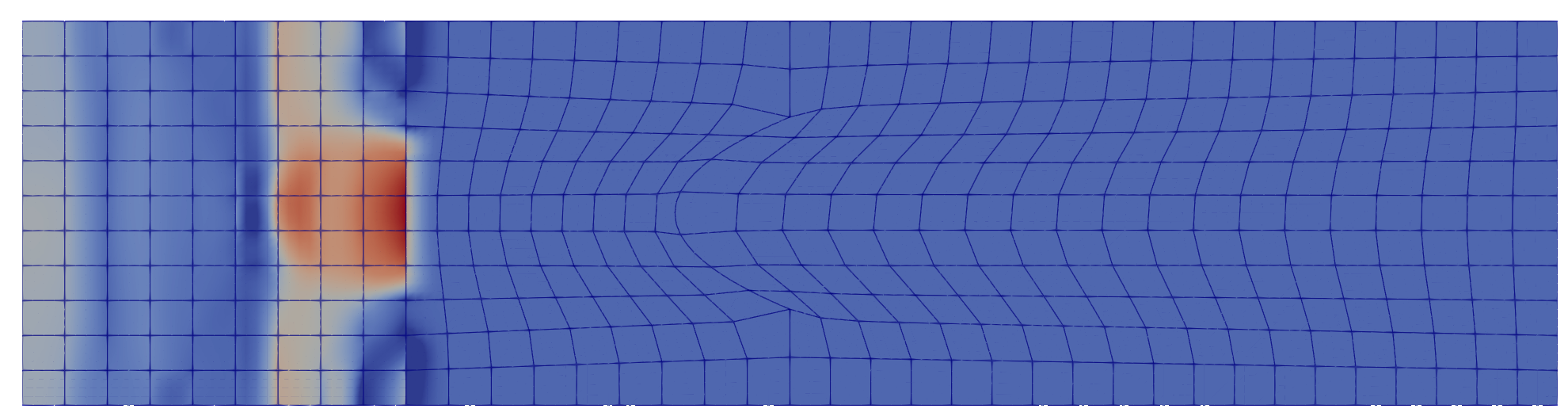}
		\end{subfigure}
		};	
		\draw(0, -2.2) node[inner sep=0, anchor=north]{\begin{subfigure}[b]{.33\textwidth}
			\centering
			\includegraphics[width=\textwidth]{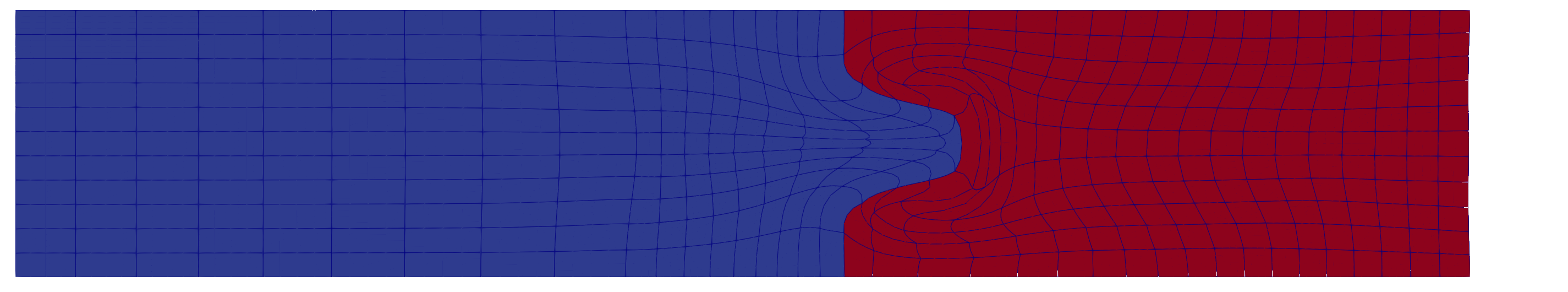}
		\end{subfigure}
		};	
		\draw(5, -2.2) node[inner sep=0, anchor=north]{\begin{subfigure}[b]{.33\textwidth}
			\centering
			\includegraphics[width=\textwidth]{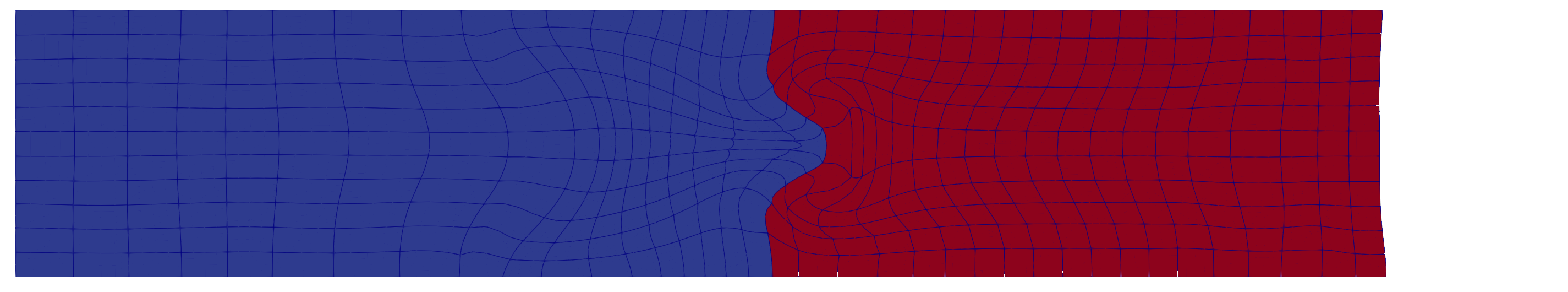}
			\caption{}
		\end{subfigure}
		};	
		\draw(10, -2.2) node[inner sep=0, anchor=north]{\begin{subfigure}[b]{.33\textwidth}
			\centering
			\includegraphics[width=\textwidth]{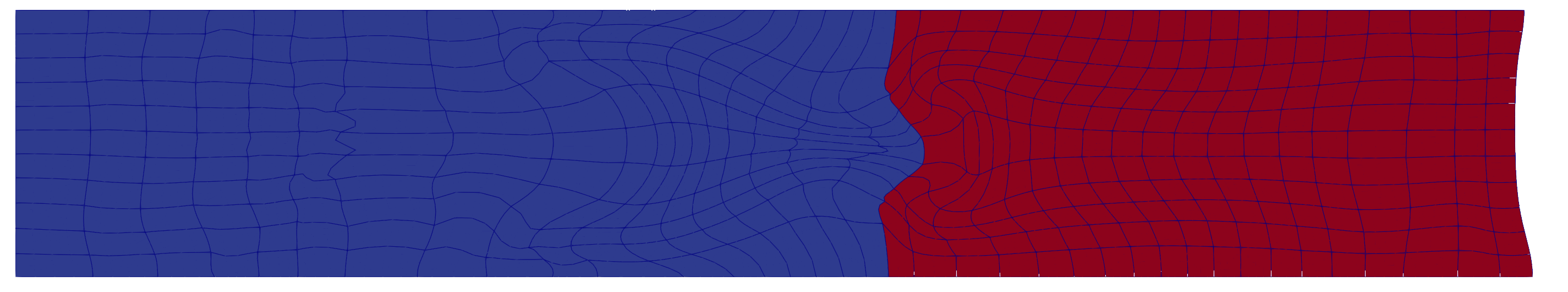}
		\end{subfigure}
		};	

		\draw(0, -4.3) node[inner sep=0, anchor=north]{
			\begin{subfigure}[b]{.33\textwidth}
				\centering
				\includegraphics[width=\textwidth]{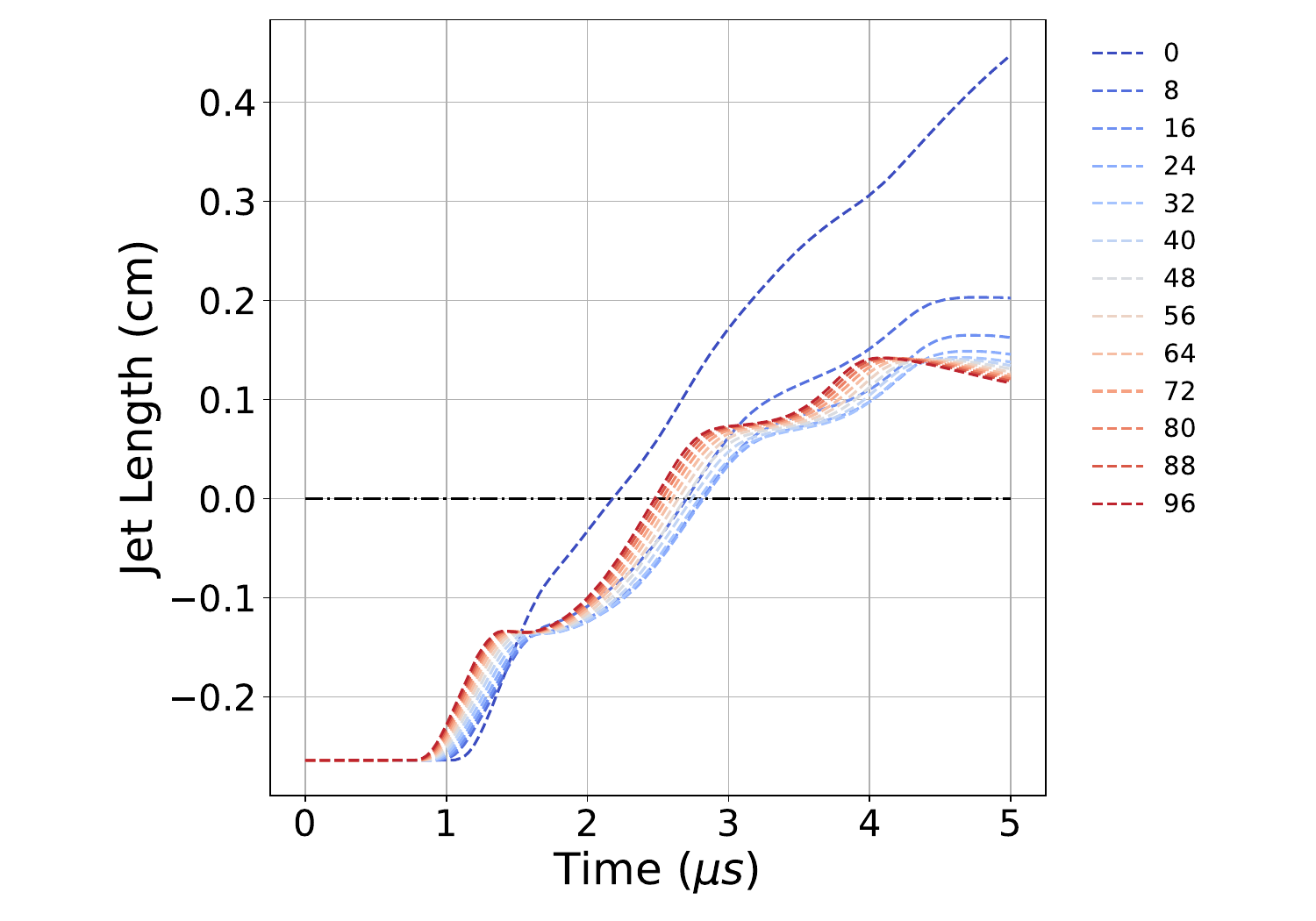}
				\caption{}
			\end{subfigure}
		};
		\draw(5.1, -4.3) node[inner sep=0,anchor=north]{
			\begin{subfigure}[b]{.33\textwidth}
				\centering
				\includegraphics[width=\textwidth]{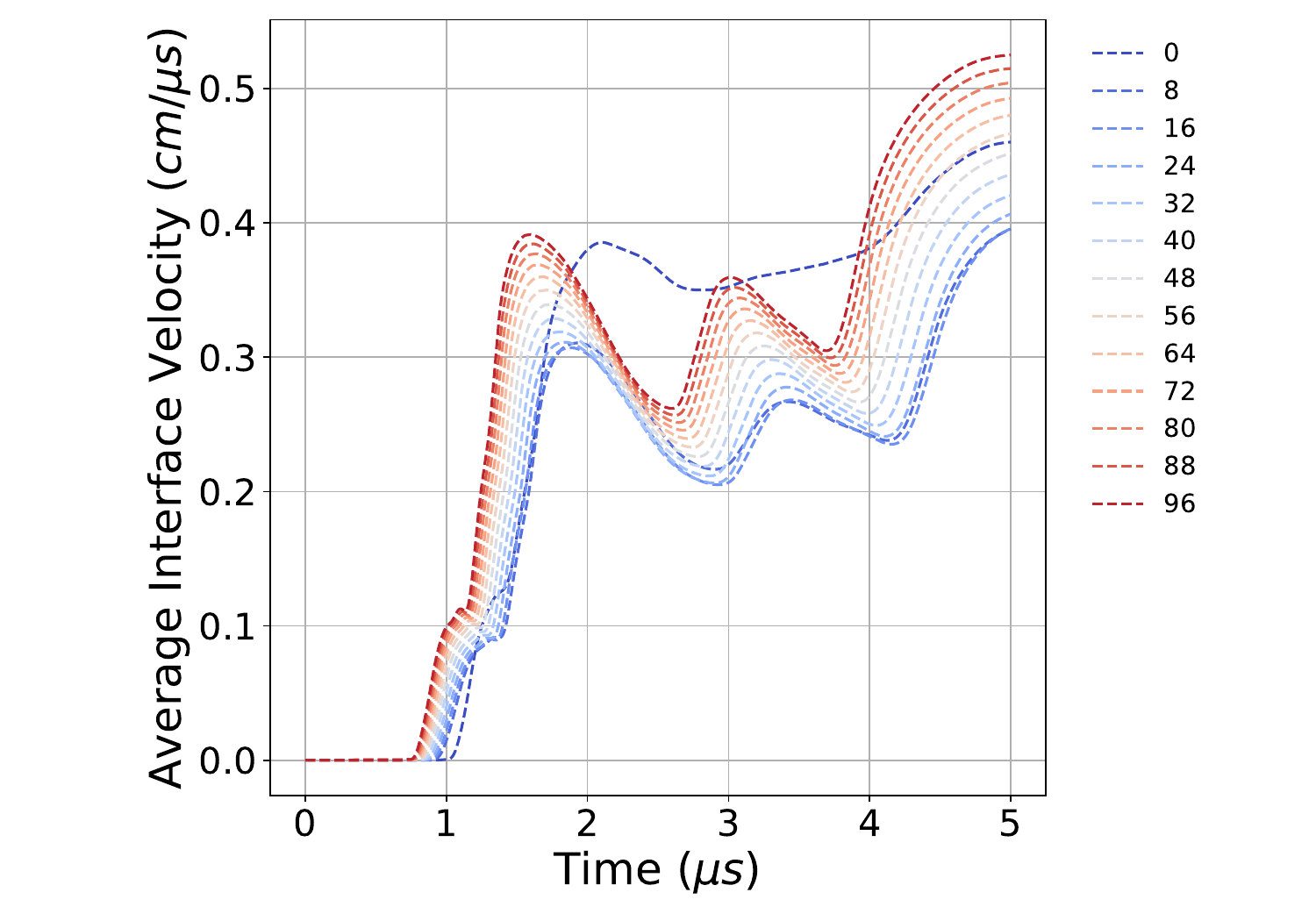}
				\caption{}
				\label{fig:interface_velocity}
			\end{subfigure}
		};
		\draw(10, -4.3) node[inner sep=0,anchor=north]{
			\begin{subfigure}[b]{.26\textwidth}
				\centering
				\includegraphics[width=\textwidth]{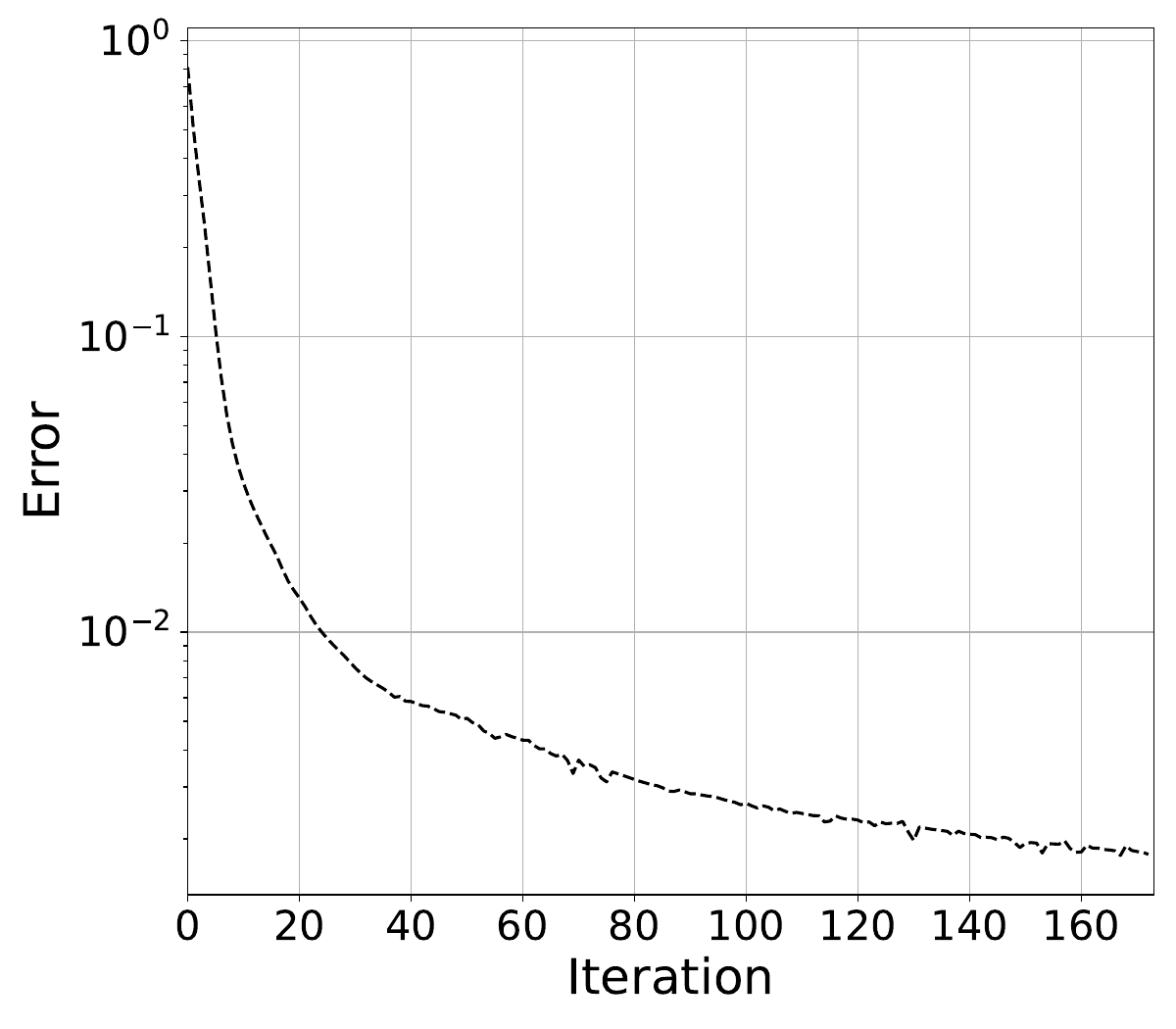}
				\caption{}
			\end{subfigure}
		};
		% Additional Writings
		\draw(-2.2, 0.3) node[anchor=north west, inner sep=0]{\footnotesize Initial Energy Field};
		\draw(-2.2, -1.9) node[anchor=north west, inner sep=0]{\footnotesize Final Deformation Profile};
		\draw [dashed, blue!70] (0.18, -3.6) -- (1.2, -4.6);
		\draw [dashed, blue!50] (5.12, -3.6) -- (1.42, -4.95);
		\draw [dashed, red] (10.12, -3.6) -- (1.45, -5.1);
	\end{tikzpicture}
	\caption{Demonstration of mitigation of RMI through the gradient descent procedure. (a) The initial energy profile at three different stages (initial guess, 8 steps, 88 steps). (b) The deformation profile at the final timestep for each of those same three stages. (c) The evolution of the jet length over time for different iterations of gradient descent. (d) The average interface velocity over time for different iterations. (e) The change in objective function for each iteration.}\label{fig:rmi_results}

\end{figure}

Figure~(\ref{fig:rmi_results}) summarizes the results of the optimization procedure described in the previous sections.

We note a few features of the solution. 
First, consider the bifurcation of the initial energy into left and right ``hot spots''. 
These regions will cause expansion in both of these locations, causing a complex shock front. 
Notably, it actually splits the shock into two different features. 
First, it will split the shock so it hits the bottom and top faces harder than the middle inclusion. 
This will temper RMI growth to flatten the interface. 
Next, a second, stronger shock wave hits the (now flattened) interface accelerating it more. 
This type of energy distribution is intuitive in hindsight, but we find it remarkable that this was discovered automatically through gradient based optimization.

The optimization algorithm can be seen to go through multiple stages. 
This is due to the competing objectives given in Equation~(\ref{eq:rmi_objective}).
Initially, the priority is interface flattening. 
Then, when that error is sufficiently low, the interface velocity is slowly increased. 
This behavior can be seen in Figure~(\ref{fig:interface_velocity}) as the final interface velocity initial is lowered, then, while maintaining a flat interface, the interface velocity is raised. 
Eventually, the interface velocity even surpasses that of the initial guess. 
Remarkably, this demonstrates that utilizing gradient optimization, RMI suppression is readily achievable without a decreasing the intensity of the drive. 
Moreover, as RMI is a general, challenging prototype for hydrodynamic behavior, this example demonstrates the ability of this class of optimization methods to be used in practical applications involving Lagrangian hydrodynamics.

\subsection{Energy Constrained Optimization}
A potential issue with the above results (with regard to the energy field) is that optimized solution may be infeasible for multiple reasons.
Two factors we consider here are (1) we want to operate within a particular energy budget and (2) the initial energies cannot be negative locally. 
The first condition follows from the fact that we do not want to produce solutions which require arbitrary large or small amounts of energy in order to produce. 
In practical problems, such as laser drives, this is represented by a maximal laser intensity allowed in the specifications. 
The second condition is a consequence of physical principals of thermodynamics.
Specifically, it is impossible to create a system with negative absolute temperature anywhere.
Thus, we want all our designs to remain in the strictly positive regime.
In mathematical terms, we have two separate conditions we would like to satisfy:
\begin{align}
	\int_{\Omega|_{t=0}} e d\Omega = C \, , \\	
	e(X, t = 0) \geq 0 \, .
\end{align}
The first enforces that the total energy remains constant in the optimization procedure and the second that the energy remains in the feasible regime.
This is a limitation, not of the optimization procedure, but of the formulation of the physics we are trying to model. 
As a result, incorporating both of these constraints into an optimization procedure is necessary to produce both feasible and practical results.
\subsubsection{Equality constraint}
Consider a constraint given by the equation
\begin{equation}
	g(y) = 0\, .
\end{equation}
We require that a perturbation $\delta y$ also satisfies the constraint, i.e.,
\[
	g(y + \delta y) = 0\, .
\]
Assuming small perturbations, we can Taylor expand the above equation to get
\begin{equation}
	\dfrac{\partial g}{\partial y} \cdot \delta y  = g_y \cdot \delta y= 0 \, .	
\end{equation}
As a result, we can define a new perturbation $\delta \tilde{y}$ such that
\begin{equation}
	\delta \tilde{y} = \left( I - \frac{1}{\left| g_y \right|^2} g_y \otimes g_y\right)	\delta y \, .
\end{equation}
For the specific case where we want to keep the total internal energy constant, we have the constraint
\begin{align}
	g(\mathbf{e}) = \int_\Omega e d\Omega 	
\end{align}
Taking the derivative as above, we have
\begin{align}
	\dfrac{\partial g}{\partial e_i} = \int_\Omega \dfrac{\partial e}{\partial e_i} d\Omega = \int_\Omega \phi_i d\Omega\, .
\end{align}
Conveniently, this derivative is constant with respect to the state variables, owing to the linearity of the energy conservation constraint. 
As a result, the application of the projection method described above will satisfy the constraint exactly.	
For shorthand, we will refer to the constraint projection as $\delta \hat{e} = P(\delta e)$.
\subsubsection{Inequality Constraint}
The second constraint we want to satisfy is that the initial energy is locally non-negative.
Because we are considering functions $e \in L^2(\Omega, t)$, we can take advantage of the property that values of $e$ within the element will be extremal at node values. 
Therefore, it is sufficient to apply the constraint on the discretization $e_i$.
The rectification of values is not a unique process; however, we choose the simplest as 
\begin{align}
	\delta \hat{e} = \text{ReLU}(e + \delta e) - e\, ,	
\end{align}
where $\text{ReLU}(x) = \max(0, x)$.
This transformation takes a perturbation $\delta e$ and energy state $e$ as inputs, rectifies the sum of the two, then returns the perturbation such that all values of the resulting energy are non-negative. 
For shorthand, we will refer to this transformation as
\begin{align}
	\delta \hat{e} = R(\delta e, e)\, ,
\end{align}

\subsubsection{Combination}\label{sec:Energy_Conservation_Combination}
There are many ways to satisfy the above two constraints in an optimization paradigm.
The standard is to define Lagrange multipliers for both constraints, then evaluate the KKT conditions to ensure feasibility of the perturbations.
We take a different approach.
Because the dimensionality of the inequality constraint can be quite high, we choose to modify the perturbations to ensure that the constraints remain satisfied.
We use the following algorithm to enforce the constraint
\begin{enumerate}
	\item Solve the adjoint problem to obtain an initial $\delta e$
	\item Set a \textit{tolerance} value for constraint violation
	\item While \textit{error} is greater than \textit{tolerance}
		\begin{enumerate}
			\item Solve $\delta e \leftarrow R(\delta e, e)$	
			\item Solve $\delta e \leftarrow P(\delta e)$
			\item Evaluate infeasibility \textit{error} $ = E(e + \delta e)$
		\end{enumerate}
\end{enumerate}
where we use an infeasibility error $E(x) = \text{max}(\text{ReLU}(-x))$ which is effectively a $L^\infty$ norm on the error. Additionally, because the projection operator $P$ is always applied second, we ensure that the resulting perturbation exactly satisfies energy conservation. 
\subsubsection{Results}\label{sec:RMI_Optimization_Constrained}
We repeat the optimization procedure described in Section~\ref{sec:RMI_Optimization} with the additional energy constraints described in Section~\ref{sec:Energy_Conservation_Combination}.

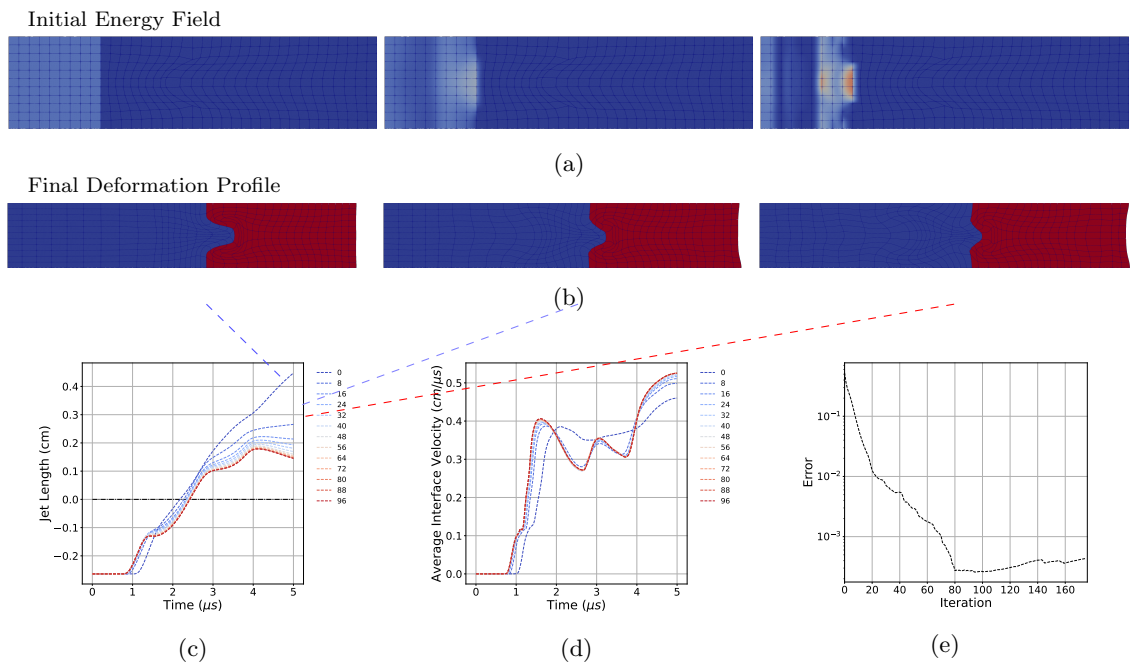
\begin{figure}[!hbt]
% \begin{figure}[H]
	\centering
	\begin{tikzpicture}
		\draw(0, 0) node[inner sep=0, anchor=north]{\begin{subfigure}[b]{.33\textwidth}
			\centering
			\includegraphics[width=\textwidth]{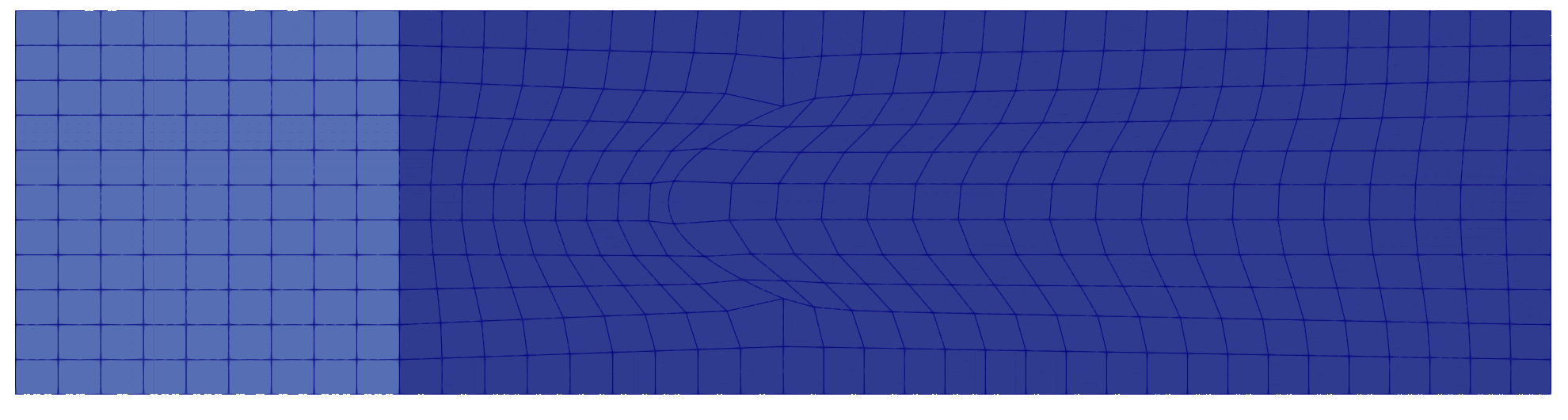}
		\end{subfigure}
		};	
		\draw(5, 0) node[inner sep=0, anchor=north]{\begin{subfigure}[b]{.33\textwidth}
			\centering
			\includegraphics[width=\textwidth]{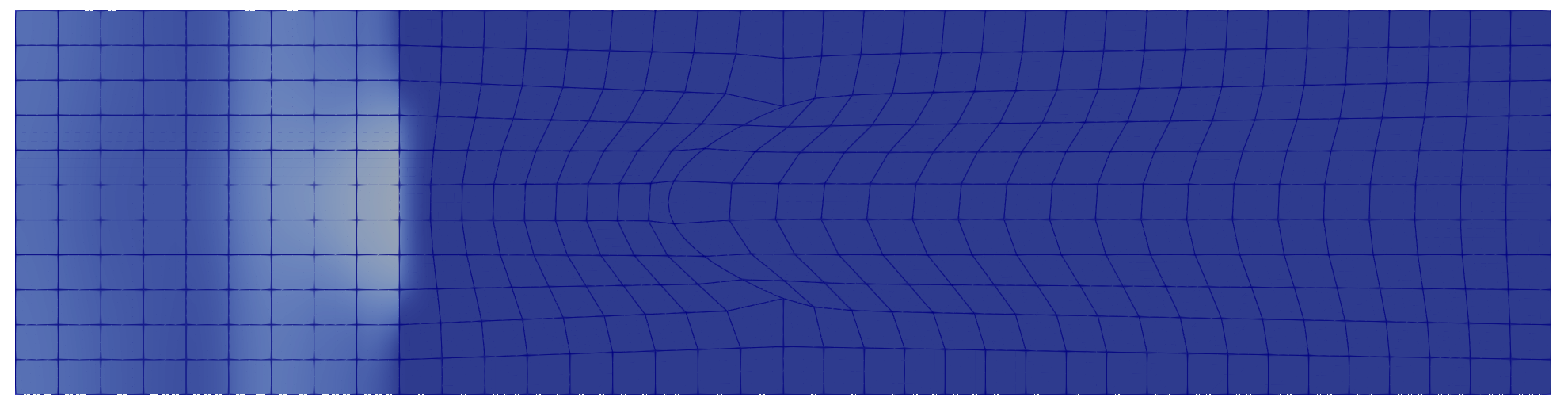}
			\caption{}
		\end{subfigure}
		};	
		\draw(10, 0) node[inner sep=0, anchor=north]{\begin{subfigure}[b]{.33\textwidth}
			\centering
			\includegraphics[width=\textwidth]{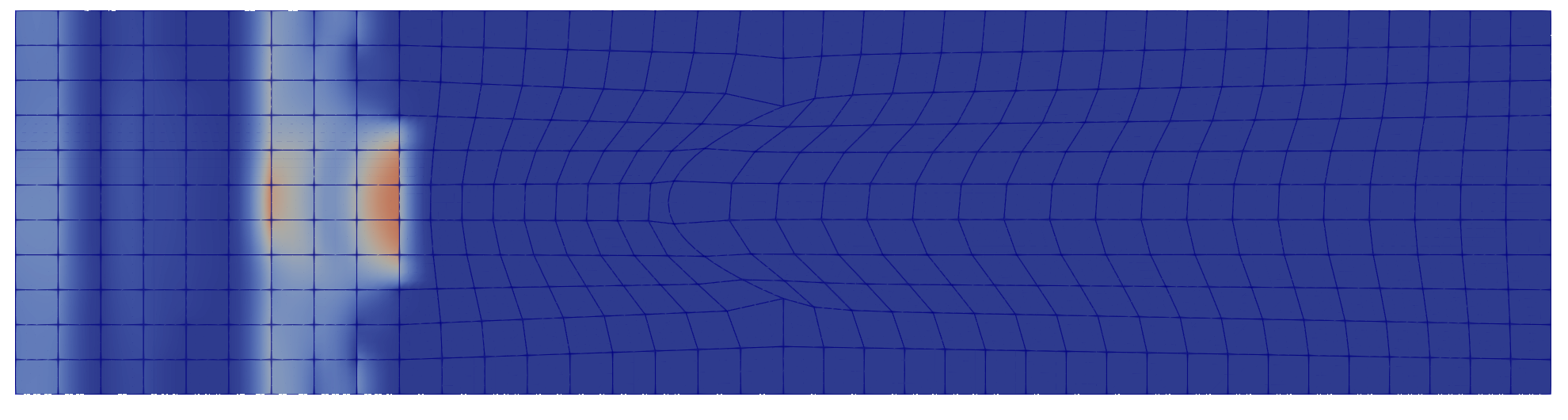}
		\end{subfigure}
		};	
		\draw(0, -2.2) node[inner sep=0, anchor=north]{\begin{subfigure}[b]{.33\textwidth}
			\centering
			\includegraphics[width=\textwidth]{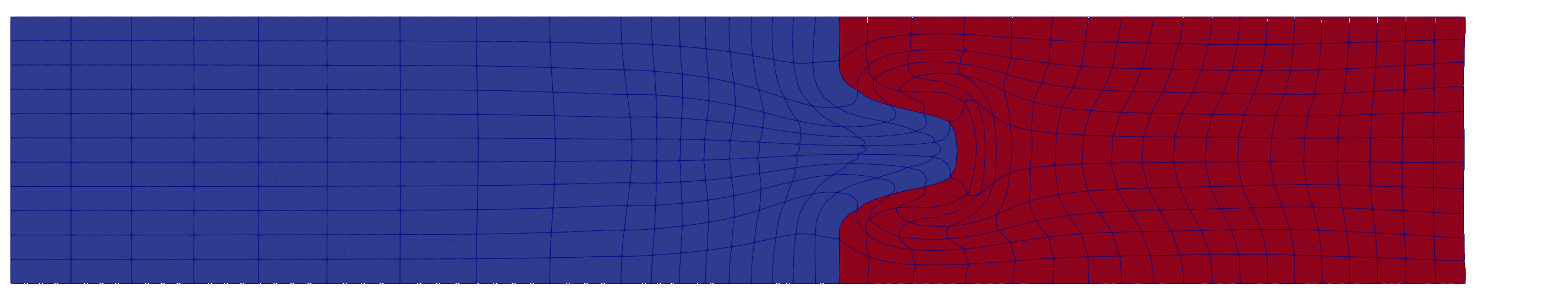}
		\end{subfigure}
		};	
		\draw(5, -2.2) node[inner sep=0, anchor=north]{\begin{subfigure}[b]{.33\textwidth}
			\centering
			\includegraphics[width=\textwidth]{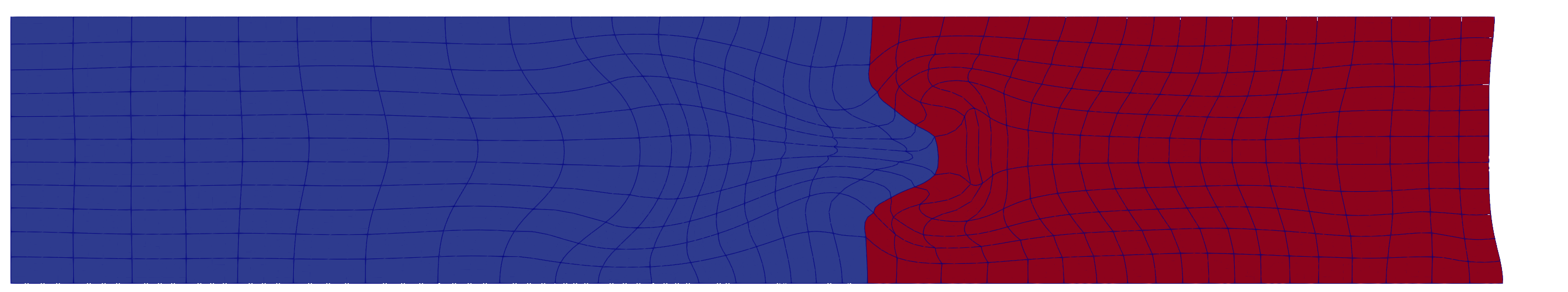}
			\caption{}
		\end{subfigure}
		};	
		\draw(10, -2.2) node[inner sep=0, anchor=north]{\begin{subfigure}[b]{.33\textwidth}
			\centering
			\includegraphics[width=\textwidth]{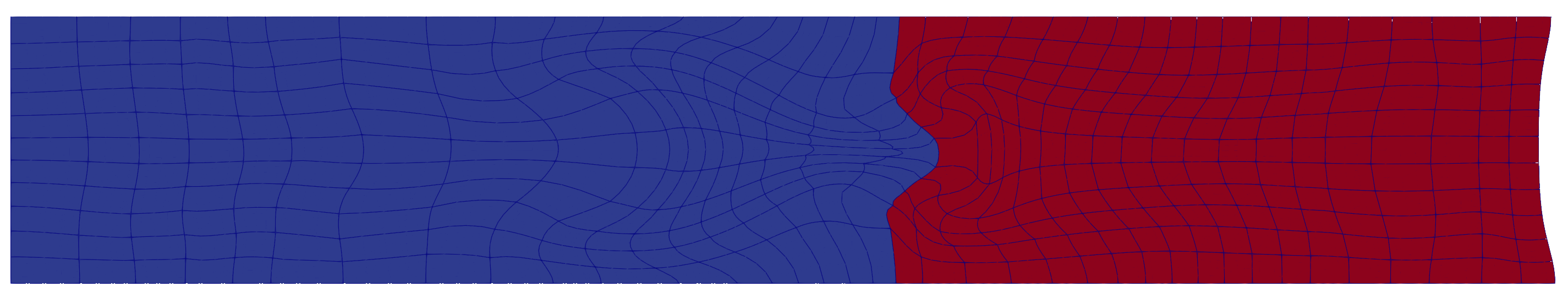}
		\end{subfigure}
		};	

		\draw(0, -4.3) node[inner sep=0, anchor=north]{
			\begin{subfigure}[b]{.33\textwidth}
				\centering
				\includegraphics[width=\textwidth]{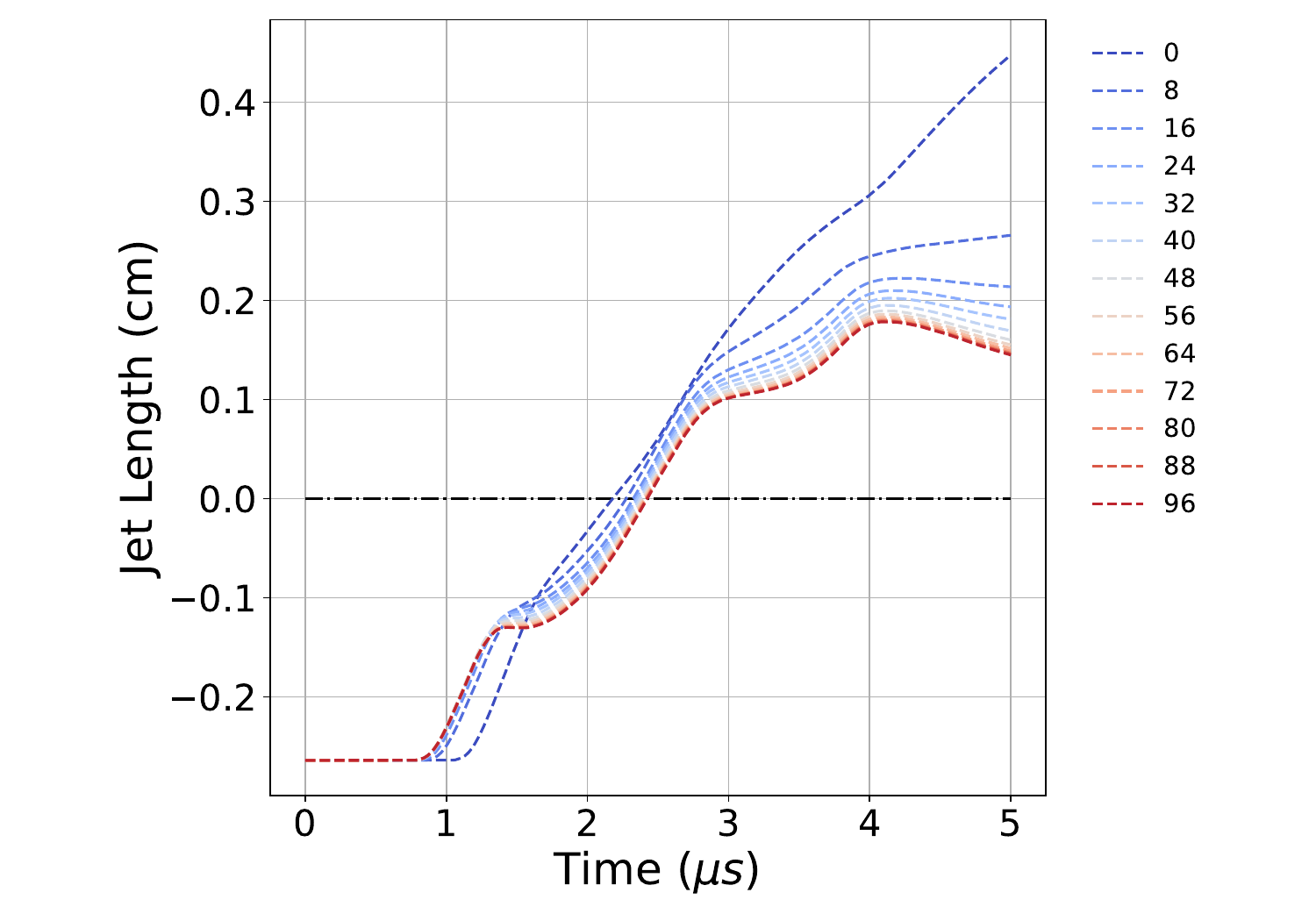}
				\caption{}
			\end{subfigure}
		};
		\draw(5.1, -4.3) node[inner sep=0,anchor=north]{
			\begin{subfigure}[b]{.33\textwidth}
				\centering
				\includegraphics[width=\textwidth]{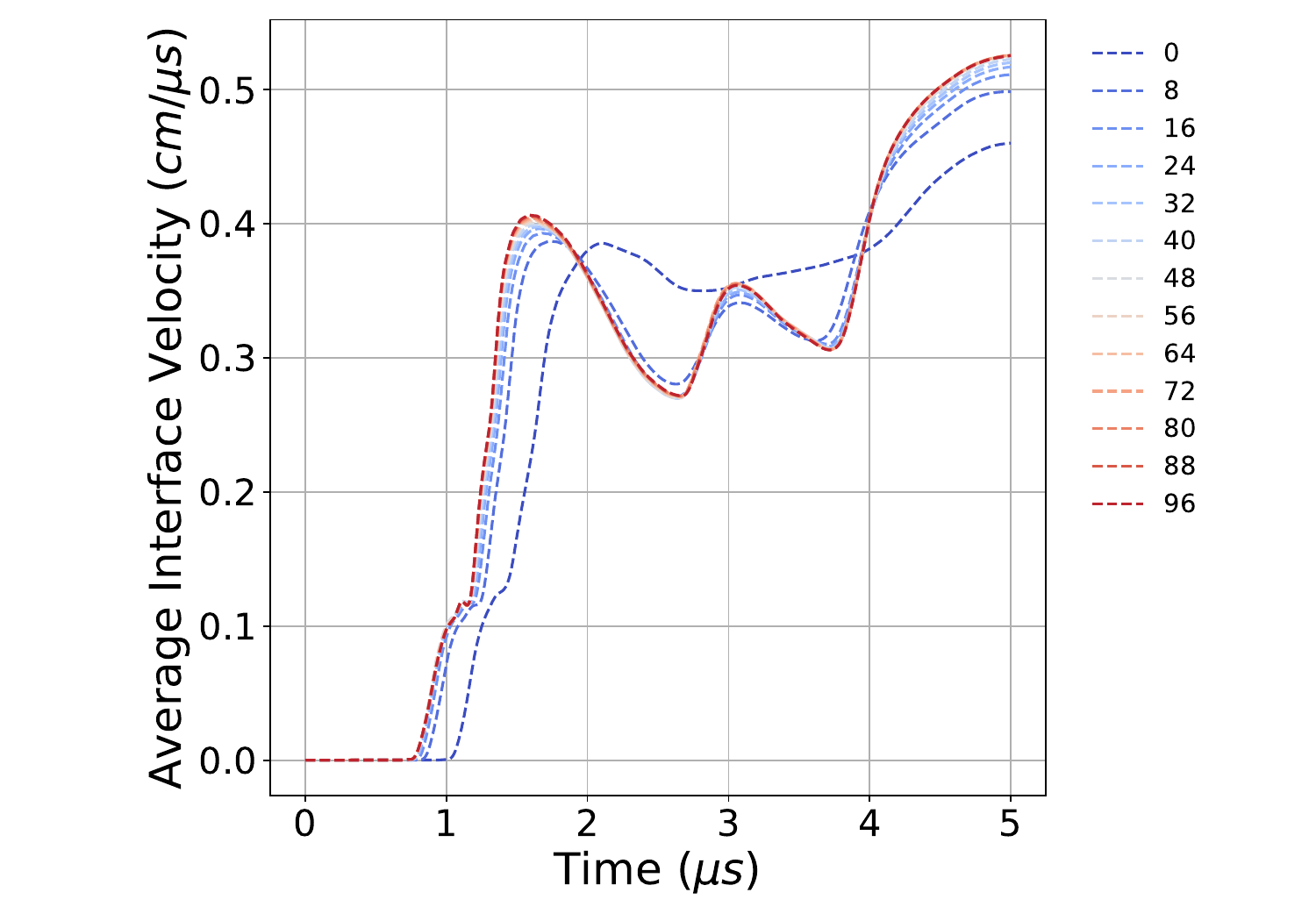}
				\caption{}
				\label{fig:interface_velocity_energy_constrained}
			\end{subfigure}
		};
		\draw(10, -4.3) node[inner sep=0,anchor=north]{
			\begin{subfigure}[b]{.26\textwidth}
				\centering
				\includegraphics[width=\textwidth]{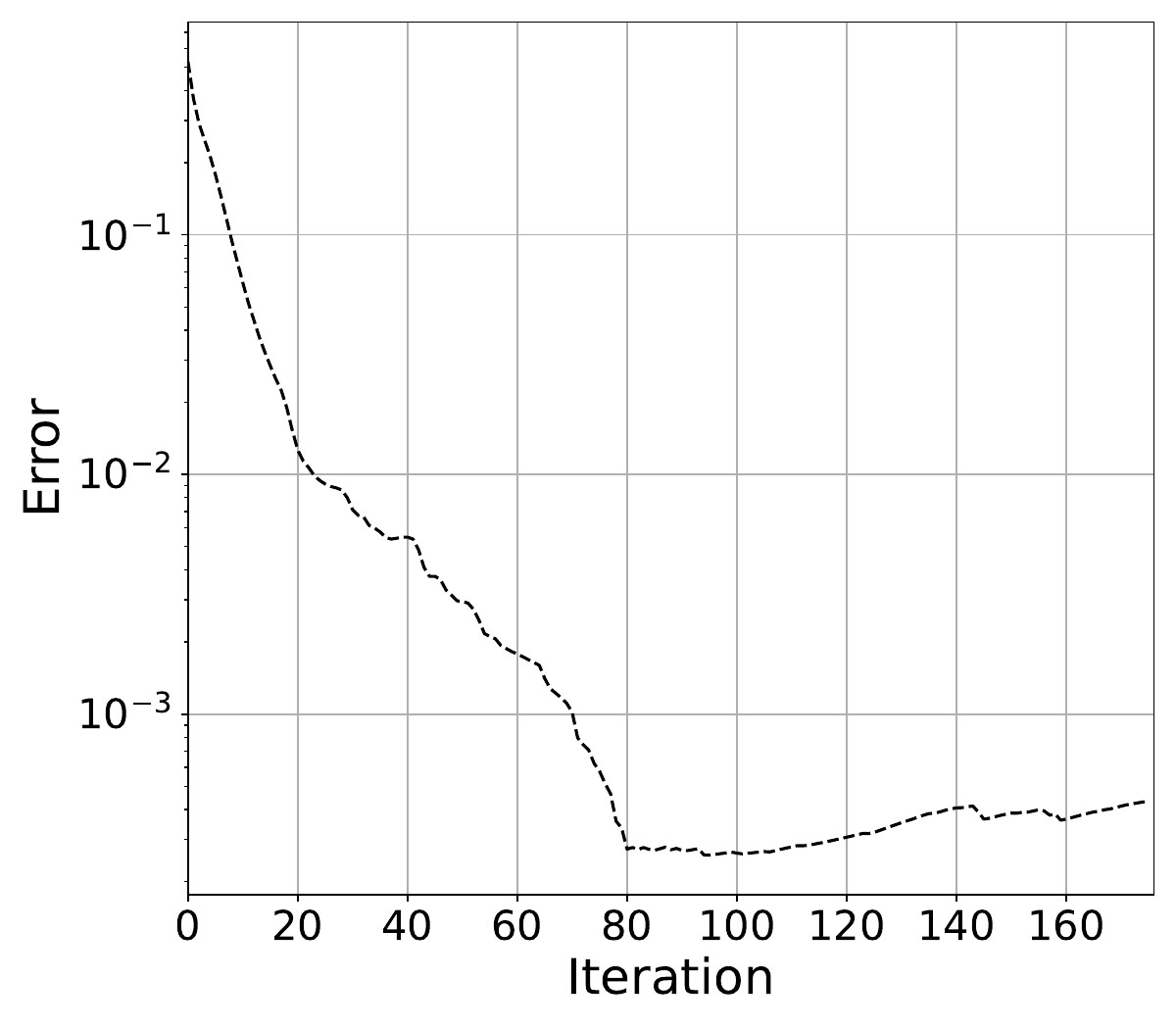}
				\caption{}
			\end{subfigure}
		};
		% Additional Writings
		\draw(-2.2, 0.3) node[anchor=north west, inner sep=0]{\footnotesize Initial Energy Field};
		\draw(-2.2, -1.9) node[anchor=north west, inner sep=0]{\footnotesize Final Deformation Profile};
		\draw [dashed, blue!70] (0.18, -3.6) -- (1.2, -4.6);
		\draw [dashed, blue!50] (5.12, -3.6) -- (1.42, -4.95);
		\draw [dashed, red] (10.12, -3.6) -- (1.45, -5.1);
	\end{tikzpicture}
	\caption{Demonstration of mitigation of RMI through the gradient descent procedure. (a) The initial energy profile at three different stages (initial guess, 8 steps, 88 steps). (b) The deformation profile at the final timestep for each of those same three stages. (c) The evolution of the jet length over time for different iterations of gradient descent. (d) The average interface velocity over time for different iterations. (e) The change in objective function for each iteration.}\label{fig:rmi_results_energy_constrained}

\end{figure}
It can be seen that the addition of energy conservation and non-negativity constraints does indeed change the structure of the optimized solution. 
In particular, the solution avoids large amounts of energy locally above and below the hot spot in the center. 
Generally, designing for experiments will often include constraints such as this in order to generate feasible designs.

\section{Conclusion}
Gradient based optimization is dramatically more efficient than gradient free optimization; the former can readily handle far greater than O(100) of degrees of freedom whereas gradient free optimization is limited to far smaller O(10 DOF) in many practical scenarios. 
In this article, we have developed a computational strategy for efficiently differentiating the predicted outputs of a high-order finite element Lagrangian hydrodynamics code with respect to its inputs via a combination of adjoint theory and automatic differentiation; further, we have applied this to a challenging problem of interface stability. This involved individually differentiating the time-stepping update, the (partially-)assembled force vector, and zero-dimensional physics (quadrature point update). 
Success in this involved many simultaneous developments (i) an efficient checkpointing algorithm and (ii) efficient vector-matrix product strategy for the time updates. (iii) We needed a proper regularization of the artificial viscosity which provides computational regularity for studying shocks -- discontinuities in the thermodynamic state of the material – which gives this area of physics challenging character. 
In order to efficiently represent the derivative of the finite element force vector, we have leveraged (iv) partial assembly to compute the action of the resulting operator; this enables efficient computation of vector-matrix products in the time-stepping and optimization. 
Finally, we have utilized an efficient implementation of automatic differentiation at the quadrature point level to account for the extensive complexity of material function calls; this allows us to use the extreme human-time efficiency in terms of functional complexity of AD while simultaneously retaining computational efficiency of adjoint methods. 

We have applied this efficient computation of gradients to optimization of the historically challenging problem in hydrodynamics of stable shock-acceleration of density interfaces; this application is critically important to major scientific challenges in fusion energy experiments such as those conducted at the national ignition facility. In these applications, the process of confining the fusion fuel to ignition conditions invariably involves launching multiple shock waves via a laser energy source which pass through density interfaces (typically diamond – DT).
We show that for a prototype problem of this phenomenon called RMI, the gradient based optimization rapidly tailors a complex spatially dependent high dimensional energy drive which simultaneously suppresses the instability and accelerates the interface to a higher velocity than the baseline case. 
In addition, we apply various constraints to our solutions with the optimization procedure in order to more accurately represent both physical and practical limitations of experimental methods. 
This simultaneous achievement is remarkable and demonstrates the value of bringing gradient based optimization to bear on problems involving computational hydrodynamics.

We close with some speculation on the future of this field; we advocate that there are many important tasks to do. 
First, we have shown this for Lagrangian Hydrodynamics; state -of-the-art codes utilize arbitrary Lagrangian-Eulerian (ALE) strategies to manage the computational mesh via the introduction of a remap step which moves the materials appropriately to a new mesh. The absence of such a strategy tends to lead to mesh tangling. A key near-term topic of research is the differentiation of this step. 
This is made challenging since formally, one must differentiate the remap, the re-mesh, and consider the multi-material case which
involves discontinuous material interfaces.
Second, hydrodynamics oftentimes does not operate alone. In general, we must track many additional state variables associated with, for instance, phase transitions, material strength, or other multiphysics. 
Tracking the dependencies in the resulting differentiation calculation in a robust and extensible manner is challenging.
Third, we believe that there are many interesting problems in hydrodynamics that should be considered, for instance, the authors intend to eventually conduct realistic simulations of  NIF laser driven shots aiming to optimize capsule shape to account for well known and persistent laser drive asymmetries. 
Finally, there is ample opportunity to apply these techniques to machine learning and uncertainty quantification.
A cheap gradient evaluation, as we have provided in this article, can be utilized by Sobolev learning strategies~\cite{czarnecki2017sobolev} and Hamiltonian Monte Carlo~\cite{hoffman2014no}.

\section{Acknowledgements}
This work performed under the auspices of the U.S. Department of Energy by Lawrence Livermore National Laboratory under Contract DE-AC52-07NA27344. This work was supported by the LLNL-LDRD Program under Project No. 21-SI-006 and Project No. 24-ERD-005. The Lawrence Livermore National Security journal number is LLNL-JRNL-872441. We thank Libby Glascoe, John Edwards, and Teresa Bailey for their outstanding programmatic leadership and advocacy for this topic. We thank Sam Mish and Robert Carson for numerous interesting discussions.

\newpage
% \printbibliography
\bibliography{bibliography}{}
\bibliographystyle{plain}

\newpage
\section{Supplemental Information}

\subsection{Mesh Convergence}
We verify the convergence with respect to the mesh by visualizing the adjoint field.
\begin{figure}[H]
	\centering
	\includegraphics[width=0.7\textwidth]{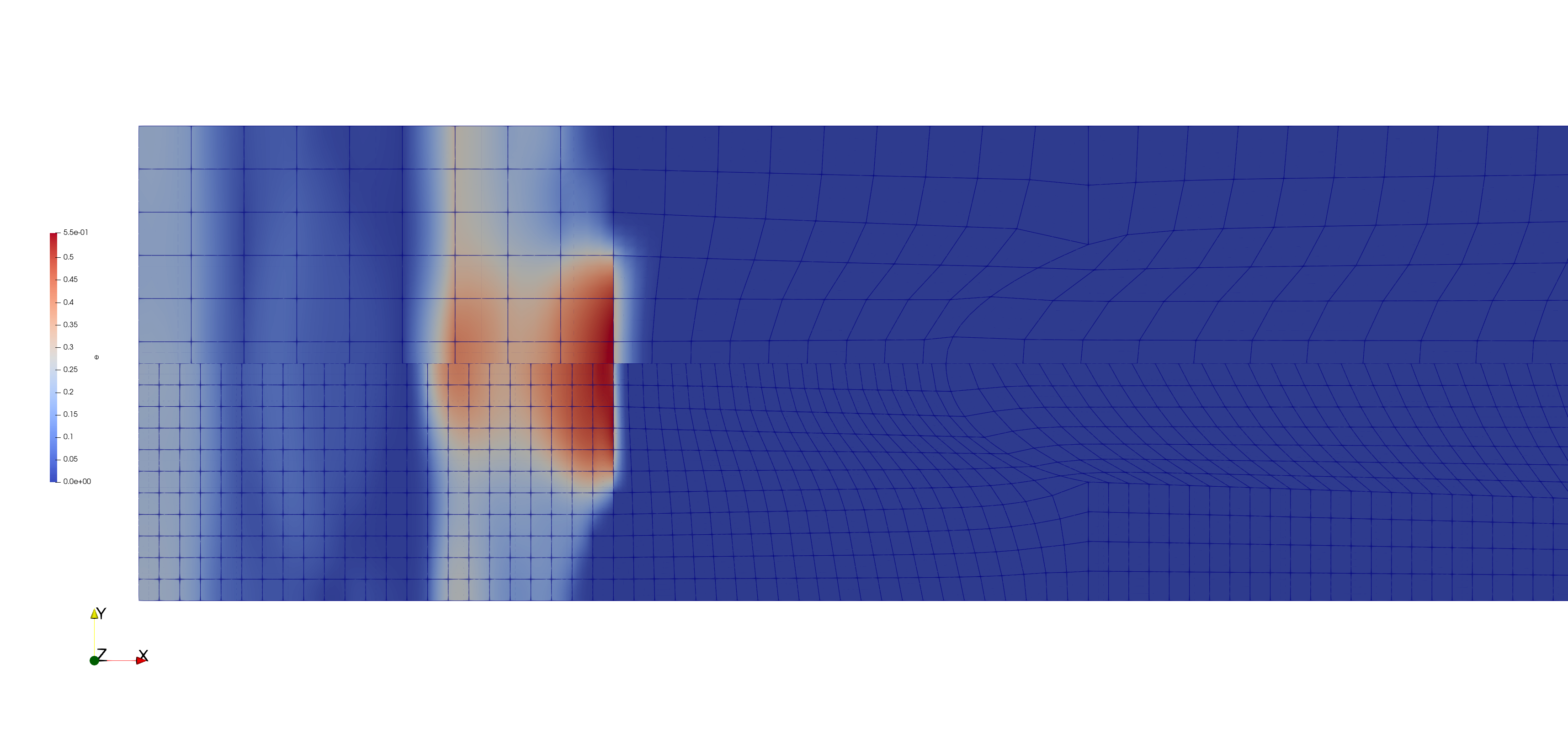}
	\caption{Solution comparison for different meshes.}\label{fig:mesh_comparison}
\end{figure}
\subsection{Lower Energy Results}
We also run the results with the initial energy initialized to $e(\Omega_1, 0) = 0.5$.
\begin{figure}[!hbt]
% \begin{figure}[H]
	\centering
	\begin{tikzpicture}
		\draw(0, 0) node[inner sep=0, anchor=north]{\begin{subfigure}[b]{.33\textwidth}
			\centering
			\includegraphics[width=\textwidth]{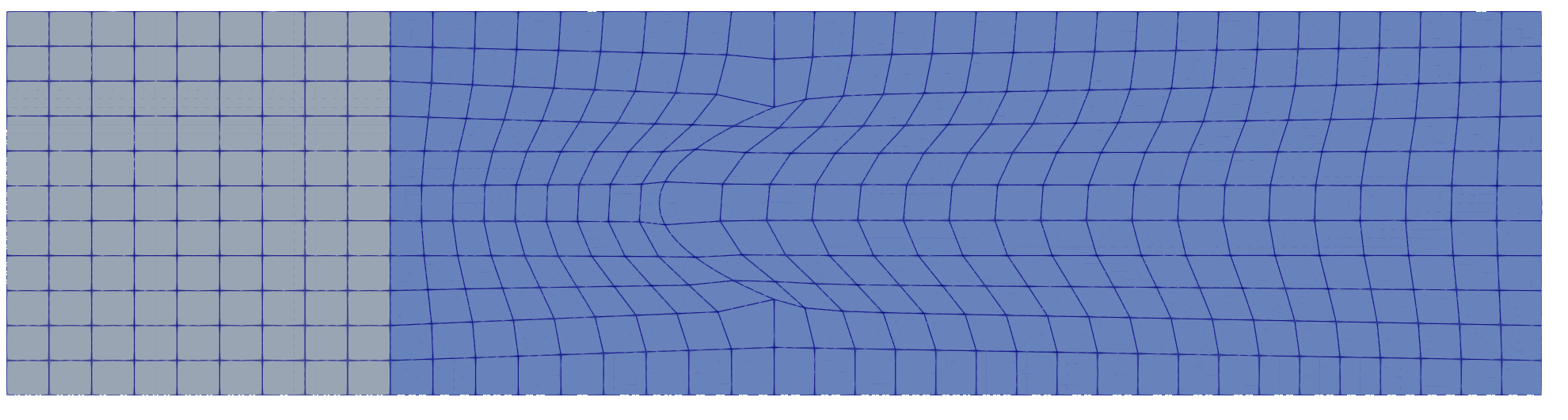}
		\end{subfigure}
		};	
		\draw(5, 0) node[inner sep=0, anchor=north]{\begin{subfigure}[b]{.33\textwidth}
			\centering
			\includegraphics[width=\textwidth]{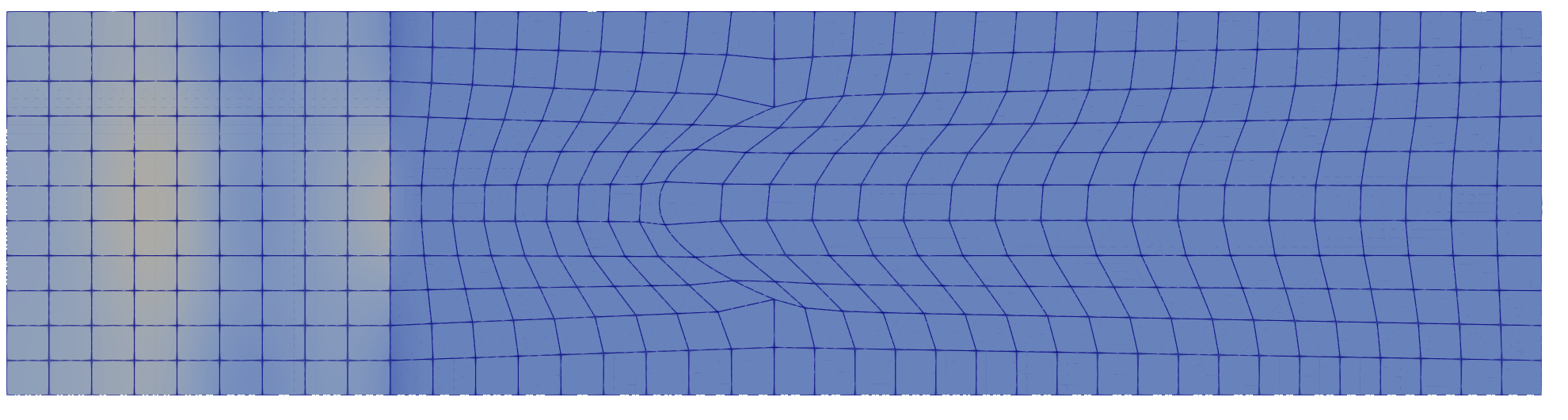}
			\caption{}
		\end{subfigure}
		};	
		\draw(10, 0) node[inner sep=0, anchor=north]{\begin{subfigure}[b]{.33\textwidth}
			\centering
			\includegraphics[width=\textwidth]{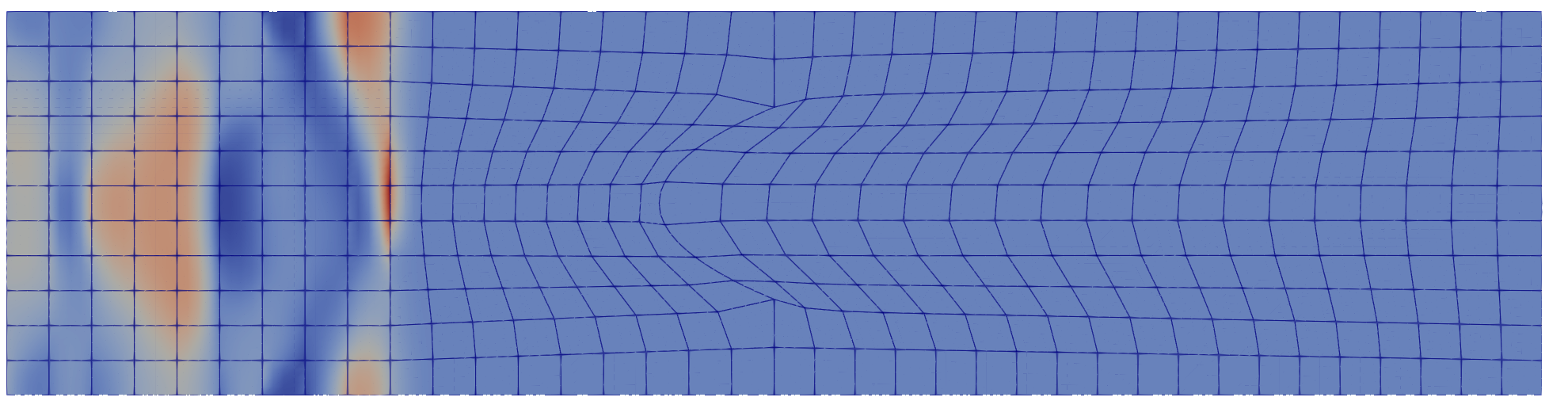}
		\end{subfigure}
		};	
		\draw(0, -2.2) node[inner sep=0, anchor=north]{\begin{subfigure}[b]{.33\textwidth}
			\centering
			\includegraphics[width=\textwidth]{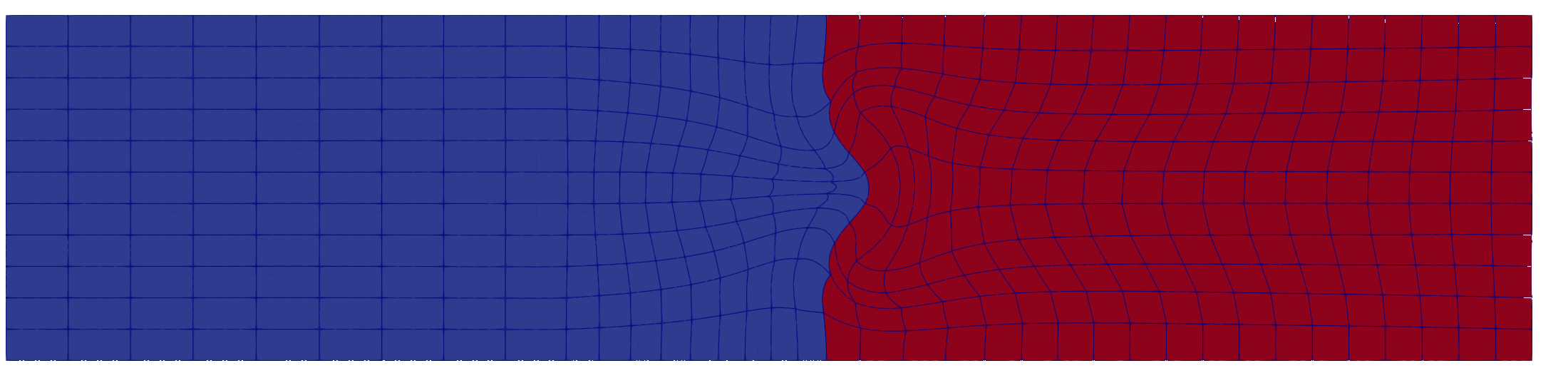}
		\end{subfigure}
		};	
		\draw(5, -2.2) node[inner sep=0, anchor=north]{\begin{subfigure}[b]{.33\textwidth}
			\centering
			\includegraphics[width=\textwidth]{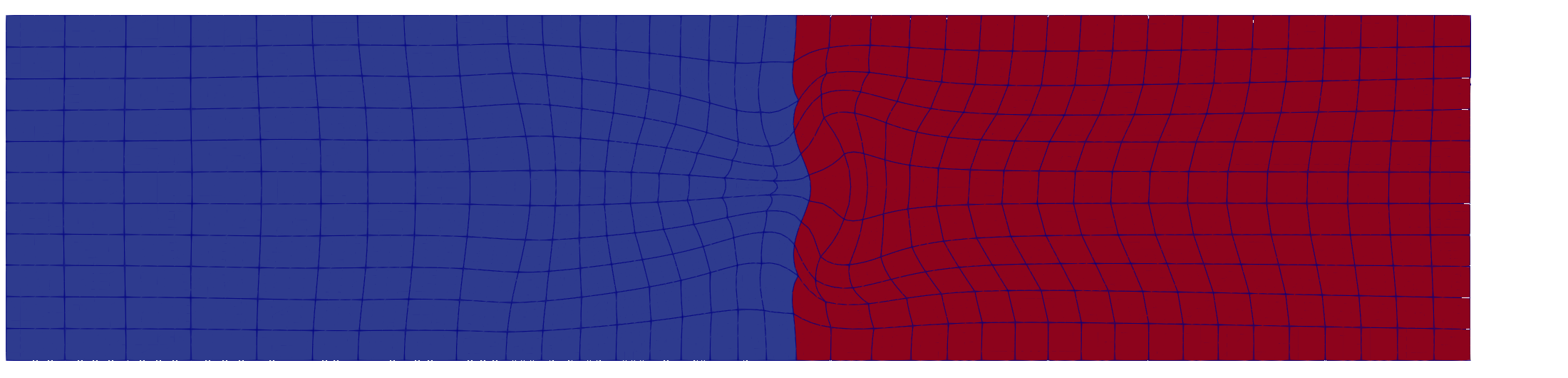}
			\caption{}
		\end{subfigure}
		};	
		\draw(10, -2.2) node[inner sep=0, anchor=north]{\begin{subfigure}[b]{.33\textwidth}
			\centering
			\includegraphics[width=\textwidth]{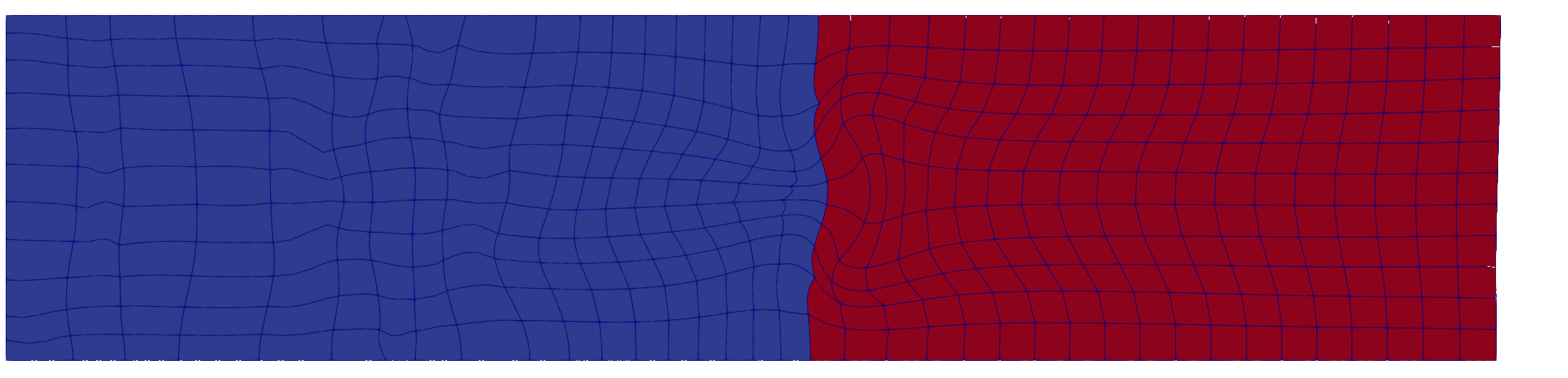}
		\end{subfigure}
		};	

		\draw(0, -4.3) node[inner sep=0, anchor=north]{
			\begin{subfigure}[b]{.33\textwidth}
				\centering
				\includegraphics[width=\textwidth]{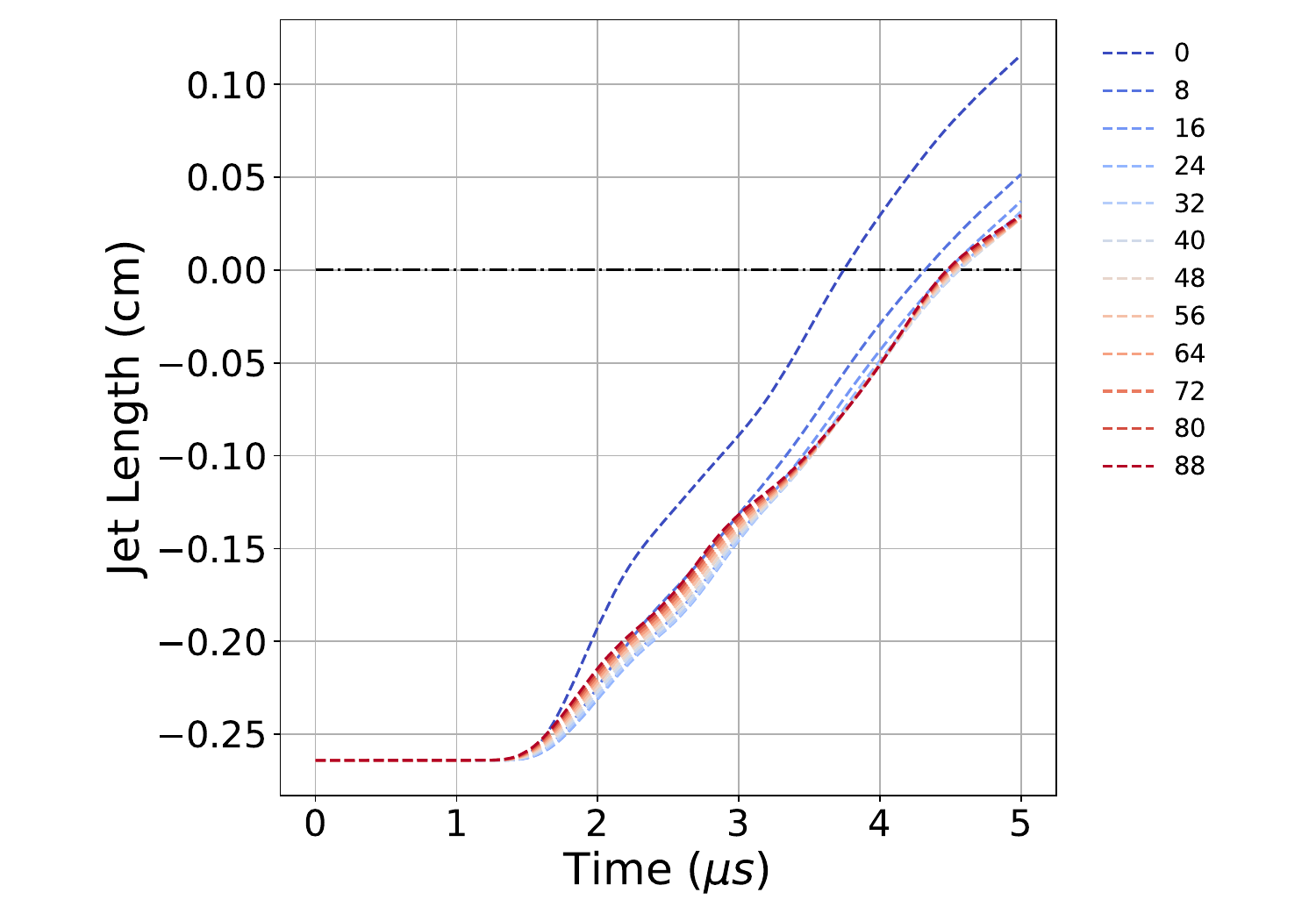}
				\caption{}
			\end{subfigure}
		};
		\draw(5.1, -4.3) node[inner sep=0,anchor=north]{
			\begin{subfigure}[b]{.33\textwidth}
				\centering
				\includegraphics[width=\textwidth]{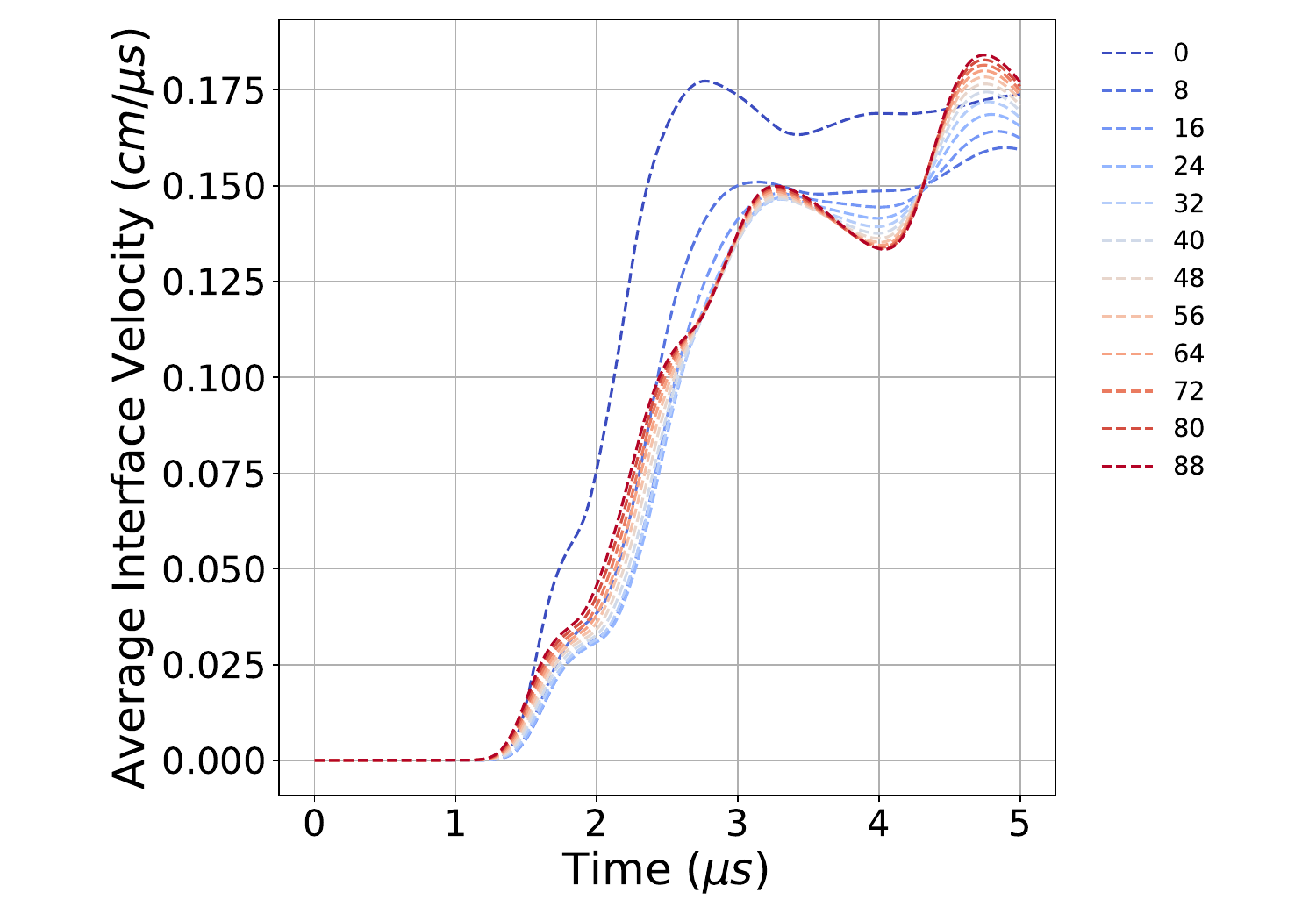}
				\caption{}
				\label{fig:interface_velocity_old}
			\end{subfigure}
		};
		\draw(10, -4.3) node[inner sep=0,anchor=north]{
			\begin{subfigure}[b]{.26\textwidth}
				\centering
				\includegraphics[width=\textwidth]{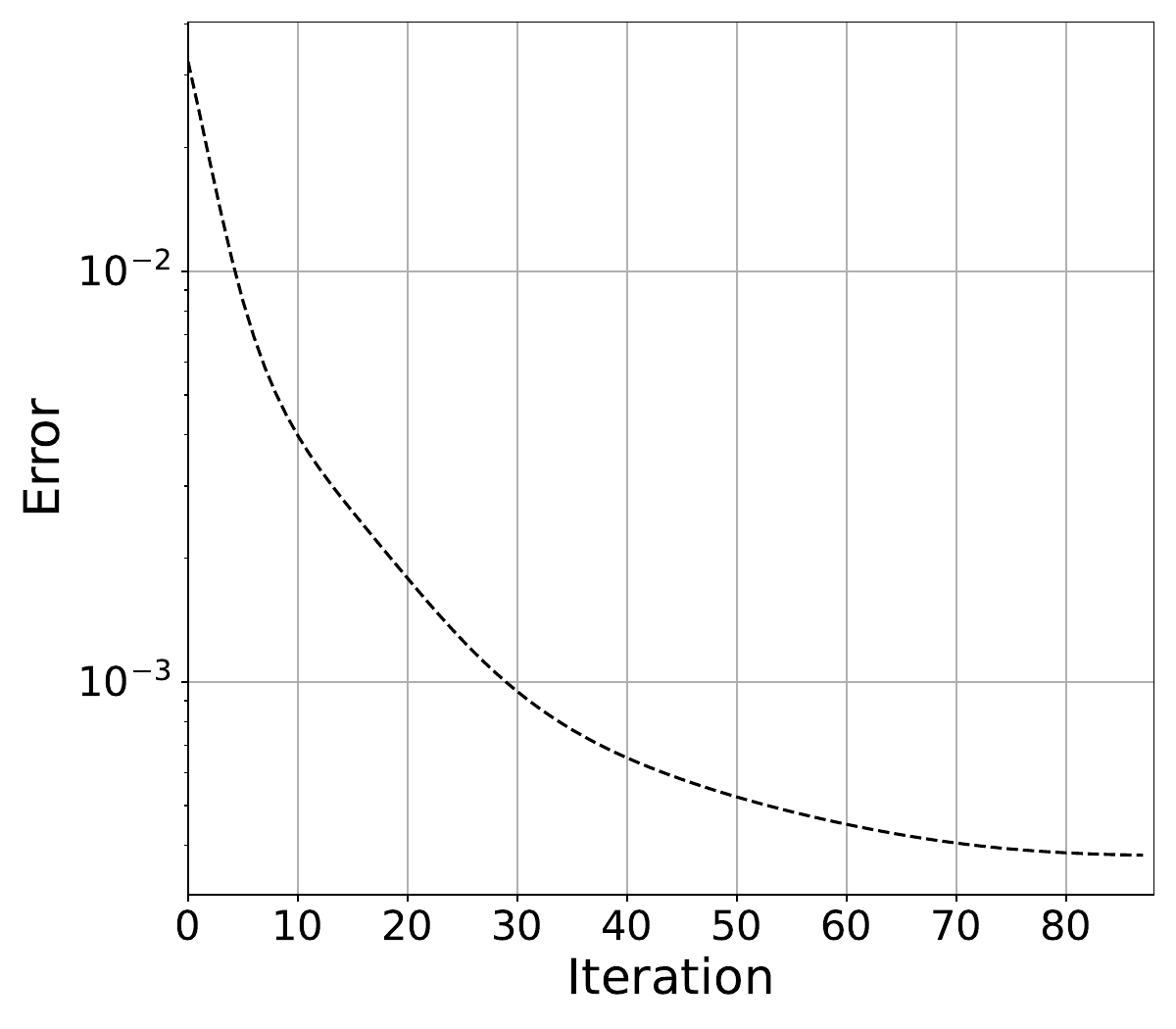}
				\caption{}
			\end{subfigure}
		};
		% Additional Writings
		\draw(-2.2, 0.3) node[anchor=north west, inner sep=0]{\footnotesize Initial Energy Field};
		\draw(-2.2, -1.9) node[anchor=north west, inner sep=0]{\footnotesize Final Deformation Profile};
		\draw [dashed, blue!70] (0.18, -3.6) -- (1.2, -4.6);
		\draw [dashed, blue!50] (5.12, -3.6) -- (1.42, -4.95);
		\draw [dashed, red] (10.12, -3.6) -- (1.45, -5.1);
	\end{tikzpicture}
	\caption{Demonstration of mitigation of RMI through the gradient descent procedure. (a) The initial energy profile at three different stages (initial guess, 8 steps, 88 steps). (b) The deformation profile at the final timestep for each of those same three stages. (c) The evolution of the jet length over time for different iterations of gradient descent. (d) The average interface velocity over time for different iterations. (e) The change in objective function for each iteration.}\label{fig:rmi_results_old}

\end{figure}
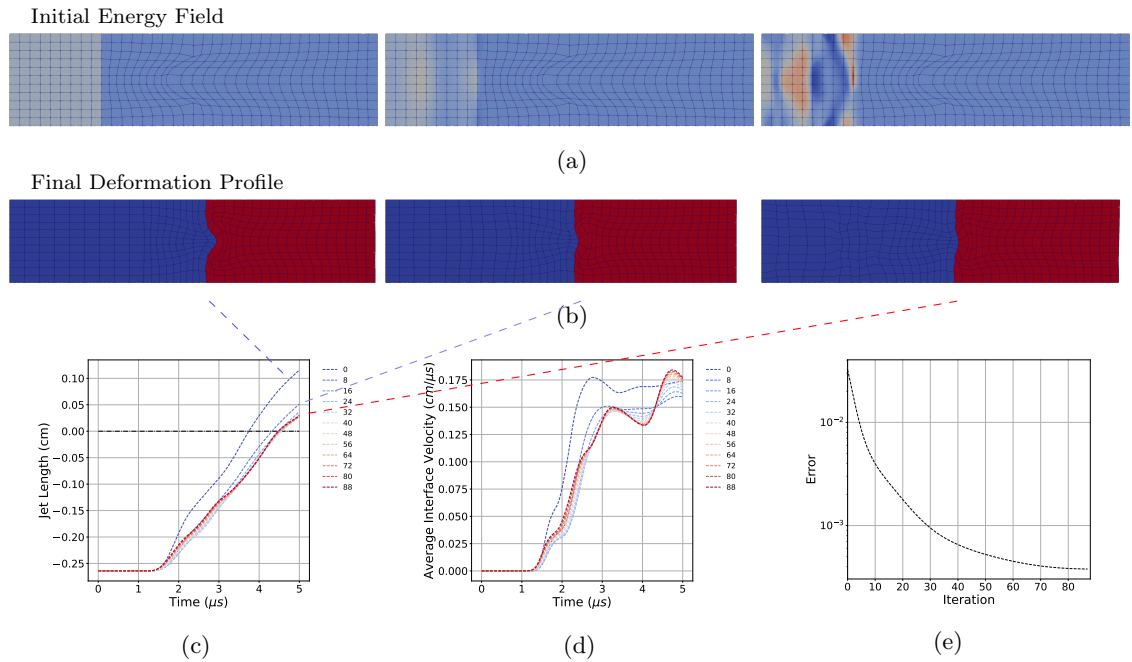

\end{document}